\documentclass[12pt]{amsart}

\usepackage[T1]{fontenc}
\usepackage{lmodern}
\usepackage{amsmath,amssymb,amsthm,mathtools}
\usepackage{abstract}
\usepackage{xcolor}
\usepackage[margin=1.08in]{geometry}
\usepackage[colorlinks=true,linkcolor=blue,citecolor=blue,urlcolor=blue]{hyperref}

\newtheorem{theor}{\hspace{1cm}{\sc Theorem}}[section]
\newtheorem{utver}[theor]{\hspace{1cm}{\sc Proposition}}
\newtheorem{sledst}[theor]{\hspace{1cm}{\sc Corollary}}
\newtheorem{lemma}[theor]{\hspace{1cm}{\sc Lemma}}
\newtheorem{conj}[theor]{\hspace{1cm}{\sc Conjecture}}
\newtheorem*{maintheor}{\hspace{1cm}{\sc Main theorem}}
\theoremstyle{definition}
\newtheorem{defin}[theor]{\hspace{1cm}{\sc Definition}}
\newtheorem{rem}[theor]{\hspace{1cm}{\sc Remark}}
\newtheorem{exa}[theor]{\hspace{1cm}{\sc Example}}

\def\C{\mathbb C}
\def\Z{\mathbb Z}
\def\Q{\mathbb Q}
\def\A{\mathbb A}
\def\Pone{\mathbb P^1}

\def\CC{({\mathbb C}^\star)}

\newcommand{\Lef}{\mathbb L}
\newcommand{\M}{\mathsf M}
\newcommand{\mot}{\mathrm{mot}}
\newcommand{\Exc}{\mathop{\rm Exc}\nolimits}
\newcommand{\ord}{\mathop{\rm ord}\nolimits}
\newcommand{\Div}{\mathop{\rm div}\nolimits}
\newcommand{\Supp}{\mathop{\rm Supp}\nolimits}
\newcommand{\Tub}{\mathop{\rm Tub}\nolimits}
\newcommand{\Res}{\mathop{\rm Res}\nolimits}
\newcommand{\lct}{\mathop{\rm lct}\nolimits}
\newcommand{\Diff}{\mathop{\rm Diff}\nolimits}
\newcommand{\id}{\mathop{\rm id}\nolimits}

\newcommand{\Var}{\mathop{\rm Var}\nolimits}

\newcommand{\im}{\mathop{\rm im}\nolimits}

\begin{document}

\begin{center}
{\Large \sc The monodromy conjecture for high-order poles, and norms of residue antiderivatives} 

\vspace{3ex}

{\sc Alexander Esterov\footnote{LIMS, aes@lims.ac.uk}}
\end{center}

\begin{abstract}
We prove the strong local monodromy conjecture for poles of order at least $n-1$ of the naive motivic zeta-function for a complex polynomial of $n$ variables. Moreover, for every such pole $s$, we construct an explicit homology $\exp(2\pi i s)$-eigencycle in a neighborhood of any given 1-dimensional exceptional stratum contributing the pole. The proof combines topological techniques developed by the author with a fact of birational geometry discovered by AI.
\end{abstract}

\tableofcontents

\section{Introduction}
The monodromy conjecture compares two kinds of invariants of singularities of a polynomial $f\in\C[x_1,\ldots,x_n]$. Every point $x\in f^{-1}(0)$ gives rise to the following objects:

-- invariants of arithmetic origin, called zeta-functions of $f$ at $x$. The simplest one is the topological zeta function (\cite{DenefLoeser1992}), an element of $\Q(s)$ determined by the Euler characteristics of the strata in the exceptional divisor of a resolution of $f$ at $x$ (see formula~\eqref{eq:topological-resolution-formula}). It refines to motivic zeta-function (\cite{DenefLoeser1998,DenefLoeser1999}), having coefficients in a localized Grothendieck ring of complex varieties. We denote the finite sets of complex poles of these zeta-functions as $P_x^{top}\subset P_x^{mot}$ (the latter's definition is not obvious, see Section~\ref{sec:poles}).

-- a topological invariant $\mathrm{Mon}_x$ -- the set of eigenvalues of the monodromy action on the ordinary homology of the Milnor fiber of $f$ at $x$ (\cite{Dimca}), and its refinement $B_x$ -- the set of roots of the Bernstein--Sato polynomial of $f$ at $x$ (\cite{MalgrangeBernstein}). We know that the union of $\mathrm{Mon}_x$ over all $x$ in a small neighborhood of 0 equals $\exp(2\pi i B_0)$.

\begin{conj}[The strong motivic local monodromy conjecture] $P_0^{mot}\subset B_0$.
\end{conj}

The conjecture is proved in the following cases.

\vspace{1ex}

i) For the pole of the maximal possible order $n$: if it exists, it is unique and equals $-\lct_0(f)$, hence is a root of the Bernstein--Sato polynomial, by \cite[Theorem~3.5 and \S3.6]{NicaiseXu}. No isolatedness assumption is needed. We retain this case below as an illustration of the eigencycle construction; Theorem~\ref{thm:bs-high-order} concerns poles of order $n-1$.

\vspace{1ex}

ii) For $n=2$, see \cite{Loeser1988,Veys1995}. 

\vspace{1ex}

iii) For several important special families of polynomials $f$, such as Newton-nondegenerate singularities satisfying Loeser's nonresonance hypotheses \cite{Loeser1990}. The weaker eigenvalue form is known for superisolated surface singularities \cite{ArtalEtAl2002} and quasi-ordinary hypersurface singularities \cite{ArtalEtAl2005}. Without Loeser's nonresonance hypotheses, Lemahieu--Van Proeyen prove the topological eigenvalue form for all Newton-nondegenerate polynomials in three variables \cite{LemahieuVanProeyen}, while Esterov--Lemahieu--Takeuchi prove it in four variables \cite{EsterovLemahieuTakeuchi,EsterovLemahieuTakeuchiCorrigendum}. Larson--Payne--Stapledon prove the local motivic eigenvalue form in arbitrary dimension for simplicial Newton polyhedra, and also for all Newton-nondegenerate polynomials in three variables \cite{LarsonPayneStapledon}. These are nearby-eigenvalue statements, rather than the Bernstein--Sato-root form above.

\begin{theor}\label{mainth}
The conjecture is valid for poles of order $n-1$ of the motivic zeta-function and its Hodge and topological specializations. 

Moreover, each term of the form $1/(s-s_0)^{n-1}$ in the partial fraction decomposition of the zeta-function comes with the coefficient made up of the contributions $C_E$ (see Section~\ref{subsec:res-motivic-bracket}, especially \eqref{eq:res-individual-RV-coefficient}) of certain 1-dimensional strata $E$ in the exceptional divisor of a resolution of $f$ at 0. The monodromy eigenvalue $\exp{2\pi i s_0}$ comes from an $(n-1)$-dimensional homology eigencycle supported near any prescribed vertical curve $E$ (i.e. projecting to 0) with $C_E\ne 0$; or, in the absence of such $E$, from an $(n-2)$-dimensional homology eigencycle supported near a point of a horizontal curve (i.e. projecting to a curve). 
\end{theor}
The "moreover" part (see the detailed version in Subsection~\ref{subsec:precise-results}) is precise as written for the topological zeta-function;  Corollary~\ref{cor:poles-individual-motivic} and Theorem~\ref{thm:poles-global-RV} give details of how it extends to the motivic zeta-function and its Hodge specialization.

\subsection{The idea of the proof} \label{ssidea} We construct explicit eigencycles, as outlined below.

1. If $f(y_1,\ldots,y_n)=y_1^{d_1}\cdots y_n^{d_n}$ is a monomial with $d={\rm gcd}(d_1,\ldots,d_n)>1$, then its Milnor fiber at 0 is a retract of the union $$C\sqcup e\cdot C\sqcup \ldots \sqcup e^{d-1}\cdot C\eqno{(*)}$$ for an order $d$ element $e$ in the complex torus $\CC^n$ and a torus coset $C\simeq\CC^{n-1}$. Note that $e$ can be chosen so that the monodromy sends each $e^k\cdot C$ to $e^{k+1}\cdot C$. Thus, for every root of unity $\zeta^d=1$, a $\zeta$-eigencycle of monodromy can be constructed as $T_\zeta:=\sum_{k=1}^d \zeta^{-k} \cdot [e^k \cdot T]$, where $T\subset C$ is a topological torus retracting $C$, and $[\ldots]$ stands for the fundamental homology cycle of a submanifold.

\begin{defin}
Apply this construction to the monomial representation of an arbitrary function $f$ near a 0-dimensional stratum $E$ of the exceptional divisor of a resolution. The resulting cycle $[T]$ is called the {\it Leray tube} of $f$ over $E$, and the eigencycle $T_\zeta$ the {\it Malgrange eigencycle} over $E$ (names motivated by \cite{Leray1959} and \cite{MalgrangeBernstein}).
\end{defin}

\begin{rem}
This construction by itself proves that $\zeta$ is a monodromy eigenvalue, unless $T_\zeta=0$ in the homology of the Milnor fiber. One possible certificate of $T_\zeta\ne 0$ is $$\int_{T_\zeta} \omega/df\ne 0 $$ for $\omega=dx_1\wedge\cdots\wedge dx_n$. This inequality in particular holds if the 0-dimensional stratum contributes the pole $-m/d$ of order $n$ (i.e. each of the exceptional components containing $E$ contributes this pole of the zeta-function). Then, for $\zeta:=\exp(2\pi i m/d)$, we have: $$\int_{T_\zeta} \omega/df=\sum_{k=1}^d \zeta^k\zeta^{-k} \int_T \omega/df=d\cdot \int_T \omega/df\ne 0.$$
Interestingly, the author did not find this proof of case (i) of the monodromy conjecture in the literature. 
\end{rem}

2. We now extend the construction of Leray tubes and Malgrange eigencycles to an arbitrary codimension $r$ open stratum $E$ of the exceptional divisor. The function $f$ restricts to a monomial on a transversal slice of $E$ at $y\in E$, giving rise to the $(r-1)$-dimensional Leray tube $T(y)$ and the Malgrange eigencycle $T_{\zeta}(y)$. Thus, in a tubular neighborhood of $E$, the Milnor fiber is a locally trivial fibration with fibers of the form $(*)$. Contracting the components of the fibers, we get a degree 
$d$ covering $$p_E:\tilde E\to E.$$ 

\begin{defin}
The {\it Leray tube} $T(\alpha)$ and the {\it Malgrange eigencycle} $T_{\zeta}(\alpha)$ over a compact oriented $p$-dimensional submanifold  $\alpha\subset\tilde E$ are $(p+r-1)$-dimensional cycles in the Milnor fiber of $f$, supported near $p_E(\alpha)$ and characterised by the property that their intersection with a transversal slice of $p_E(\alpha)$ at $y\in p_E(\alpha)$ equal $T(y)$ and $T_{\zeta}(y)$ respectively (see Section~\ref{sec:tubes}, especially Proposition~\ref{prop:general-Fourier-Leray}, for a construction of such cycles, applicable to any singular cycle $\alpha$).
\end{defin}
In particular, the construction of the Malgrange eigencycle $T_{\zeta}(\alpha)$ proves that $\zeta$ is a monodromy eigenvalue, unless $T_{\zeta}(\alpha)=0$ in the homology of the Milnor fiber. A possible certificate that $T_{\zeta}(\alpha)\ne 0$ for $p=n-r=\dim E$ is

$$\int_{T_{\zeta}(\alpha)}\omega/df\ne 0.\eqno{(**)}$$

As before, if the exceptional components containing $E$ contribute the same pole $-m/d$, then choosing $\zeta:=\exp(2\pi i m/d)$ reduces $(**)$ to 

$$\int_{T(\alpha)}\omega/df\ne 0.\eqno{(***)}$$

This integral of the Gelfand-Leray form $\omega/df$ on the Milnor fiber $f=t$, as a function of $t\to 0$, has the leading asymptotics 
$t^{m/d-1}(\int_\alpha \omega_E+o(1))$ for a certain form $\omega_E$ on $\tilde E$, called (up to normalization -- see \eqref{eq:general-tube-period} and Subsection~\ref{subsec:residue-form}) the {\it residue form}.

\begin{rem}
Assume $E$ is a curve, and the residue 1-form $\omega_E$ is not exact. Then there exists a loop $\alpha\subset\tilde E$ with $\int_{\alpha}\omega_E\ne 0$, so we have $(***)$ and $(**)$, so the Malgrange eigencycle $T_\zeta(\alpha)\ne 0$ proves that $\zeta$ is a monodromy eigenvalue.
\end{rem}

3. It remains to prove that, once the vertical curve $E$ projecting to 0 has $C_E\ne 0$ (in the notation of Theorem \ref{mainth}), the residue 1-form $\omega_E$ is never exact. Towards the contradiction, assume it is $dg$ for some function $g$ on $\tilde E$. From the equivariance of $\omega_E$ under the deck transformation of the covering $p_E:\tilde E\to E$, the values of $g$ over a point $y\in E$ form the vertices $g_1(y),\ldots,g_d(y)$ of a regular polygon with centre $c(y)$ in ${\mathbb C}$. Thus the following function $G$ on $E$ is correctly defined: 
$$G:E\to{\mathbb CP}^1,\qquad G(y):=\prod_{k=1}^d (g_k(y)-c(y)).$$
The ramification profile of this function expresses in terms of combinatorics and numerology of the resolution (Lemma~\ref{lem:res-exactness-rigidity} and formulas~\eqref{eq:res-norm-divisors}--\eqref{eq:res-norm-local-degrees}), so universal properties of ramification profiles (such as the Riemann--Hurwitz formula) impose restrictions on the resolution. In particular, there is an easy contradiction between the three restrictions on the resolution of our singularity:

1) restrictions coming from universal properties of ramification profiles applied to $G$;

2) the assumption on the contribution $C_E\ne 0$ in Theorem \ref{mainth};

3) the following general fact about resolutions.

\begin{theor}[Referred to as the corner theorem; proved as Proposition~\ref{prop:mmp-corner}]\label{thmmp}
Given a projective embedded log resolution $p$ of a singularity of the zero locus of $f:\C^n\to\C$, and a component $E_i$ of the total transform $p^{-1}\{f=0\}$, let $N_i$ and $\nu_i-1$ be the multiplicities of $E_i$ in the divisors $\Div(f\circ p)$ and $\Div(p^*\omega)$ respectively. If $E_1\cap\ldots\cap E_n$ intersect at a point $p$, and $\nu_2/N_2=\cdots=\nu_n/N_n=:q$, then $\alpha_p:=\nu_1-qN_1$ does not exceed 1.
\end{theor}

The contradiction between (1-3) essentially proves the eigencycle part of Theorem \ref{mainth} for the topological zeta-function. Refining it to the motivic and Bernstein-Sato statements follows known technologies in Sections~\ref{sec:poles} and~\ref{sec:bernstein-sato}, respectively.

\subsection{The role of AI} By August 2026, the author had all the material outlined above in Subsection \ref{ssidea}, 
except for Theorem \ref{thmmp} (this corresponds to the material up to Lemma \ref{lem:res-norm-local-degrees} inclusively in the subsequent sections of the paper), plus the following:

\vspace{1ex}

$3^{\rm top}$) certain topological considerations, which the author intended to use instead of (3) above to get contradiction with (1) and (2). (These considerations along with some of their successful applications will be published as a separate text soon.)

\vspace{1ex}

The input of (1,2,$3^{\rm top}$) was given to ChatGPT Astra 6.0 (OpenAI) with an instruction to try to synthesize them into the proof of Theorem \ref{mainth}, by finding a contradiction between the numerous and messy numerical consequences of (1,2,$3^{\rm top}$). The AI model generated a proof, but, instead of using ($3^{\rm top}$), suddenly invented (3), with a proof based on deep MMP techniques. 

This invention has background in the literature. The plane-curve inequality underlying the $n=2$ case of Theorem \ref{thmmp} is Proposition~II.3.1 of \cite{Loeser1988}; it is a well known cornerstone of the proof of the monodromy conjecture for plane curves. Its extension to arbitrary $n$ relies on the ideas of the fundamental work \cite[\S2]{NicaiseXu} about the zeta pole of order $n$. That is why the author was starting the AI experiment with the belief that such an extension was unrealistic (since otherwise the experts in MMP working on the monodromy conjecture could have discovered it long before). 

The AI's original proof was much longer (95 pages), covered only the $n=3$ isolated case, and relied on the specifics of threefolds: it is attached to the arXiv publication as a non-human-verified output \cite{ChatGPTAstra} for the record (I redacted it only to delete the  outline of ($3^{\rm top}$) in the last paragraphs of Acknowledgements and of Section 1, which the AI model included in its text to state that it did not use that input). Further substantial human work and exchange with AI generalized and simplified the text to its current form.

Besides this genuine discovery, ChatGPT was routinely used to bookkeep computationally intensive parts of the work, draft technical reasonings from condensed outlines, etc. In particular, it notified the author that the construction of Malgrange eigencycles is very close to a recent generalization \cite[Proposition~7.1 and \S7.3]{Blanco} of \cite{Loeser1988} for plane curve singularities, and the construction of the function $G$ is very close to what is called the norm (compare Definition~\ref{eq:res-norm-definition}) in the cyclic extension $\C(\widetilde C)/\C(C)$ (see e.g. \cite[Remark~5.47]{MilneFields}), though the two have probably never met before in the literature.

\subsection{Precise form of the results}\label{subsec:precise-results}

Let $f\in\C[x_1,\ldots,x_n]$ be nonconstant, with $f(0)=0$.  For
$x\in f^{-1}(0)$, write $M_x$ for the Milnor fiber of the singularity of $f$ at $x$ and $\M_x$ for
the monodromy operator on ordinary homology of $M_x$ with complex
coefficients (see e.g. \cite{Dimca} for details in our context).
We put
$$
  H_j(M_{x};\C)_\lambda
  :=
  \ker\bigl(\M_x-\lambda\id\bigr),
  \qquad \lambda=e^{-2\pi iq},\quad q\in\Q_{>0}.
$$

Fix a projective algebraic embedded resolution
$$
  \rho: X\to U\subset\A^n_{\C},
$$
where $U$ is a Zariski-open neighborhood of $0$, $X$ is smooth,
$\rho$ is an isomorphism over $U\setminus f^{-1}(0)$, and the reduced
total transform $(f\circ\rho)^{-1}(0)_{\mathrm{red}}$ has simple
normal crossings. We also require $\rho^{-1}(0)$ to be a union of some
irreducible components of this divisor. This can be arranged by first blowing up
the origin and then resolving. This additional condition makes the fiber
over $0$ a union of open resolution strata, as used in the local zeta
formula~\eqref{eq:topological-resolution-formula}. Write
\begin{equation}\label{eq:intro-data}
  \Div(f\circ\rho)=\sum_{i\in J}N_iE_i,
  \qquad
  \Div\bigl(\rho^*(dx_1\wedge\cdots\wedge dx_n)\bigr)
     =\sum_{i\in J}(\nu_i-1)E_i.
\end{equation}
Strict transforms are included; a factor occurring in $f$ with
multiplicity $m$ gives a strict-transform component with data $(N_i,\nu_i)=(m,1)$.

An irreducible component $C$ of an intersection of $n-1$ components $E_i$ with
common ratio $\nu_i/N_i=q$ is an \emph{equal-ratio curve with ratio $q$}.  We call it
\emph{vertical} if $C\subset\rho^{-1}(0)$, and \emph{clean} if no further
ratio-$q$ component meets it. Here ``vertical'' means contracted to the
origin; ``clean'' is our shorthand for the stated condition. A curve
meeting $\rho^{-1}(0)$ and mapping onto a curve is called
\emph{horizontal}. For a clean curve $C$ with ratio $q_C$, let $P_C$ be the finite
set of intersections with the remaining total-transform components and put
$$
  \alpha_p:=\nu_{j(p)}-q_CN_{j(p)}\ne0
  \qquad (p\in P_C).
$$
Its individual principal motivic contribution is
\begin{equation}\label{eq:intro-Bmot}
  B_{\mot}(C)
   := [C\setminus P_C]
      +\sum_{p\in P_C}\frac{\Lef-1}{\Lef^{\alpha_p}-1}.
\end{equation}
Here $\Lef=[\A^1]$ is the class of the affine line. Rational powers
of $\Lef$ are formal: for a common denominator $d$ of the exponents,
one adjoins an invertible symbol $\Lef^{1/d}$ satisfying
$(\Lef^{1/d})^d=\Lef$, and puts
$\Lef^{a/d}=(\Lef^{1/d})^a$.
The coefficient ring and the principal Rodrigues--Veys pole convention
are made precise in Section~\ref{sec:poles}.

Since $\rho$ is an isomorphism over $f\ne0$, we identify each Milnor
fiber $M_{x,t}$, $t\ne0$, with its inverse image in $X$. The phrase
\emph{a homology class is supported near the set $C$} means that, for every neighborhood $V$ of $C$
and all sufficiently small $t\ne0$ relative to $V$, the class has a representative singular cycle in
$V\cap\rho^{-1}(M_{0,t})$. 

\begin{maintheor}[Detailed form of Theorem~\ref{mainth}]
Let $n\ge2$ and $q\in\Q_{>0}$.
\begin{enumerate}
\item Every prescribed vertical equal-ratio curve $C$ carries a nonzero
Malgrange eigencycle in $H_{n-1}(M_0;\C)_\lambda$, represented in every
prescribed neighborhood of $C$, provided that either $C$ is not clean,
or $C$ is clean and $B_{\mot}(C)\ne0$.

\item If $-q$ is a principal Rodrigues--Veys pole of exact order $n$, then
some point over the origin lies on $n$ ratio-$q$ components.  Every such
point carries a Malgrange eigencycle in
$H_{n-1}(M_0;\C)_\lambda$, supported in an arbitrarily small
neighborhood of that point.

\item If $-q$ is a principal Rodrigues--Veys pole of exact order $n-1$,
then one of the following alternatives holds:
  \begin{enumerate}
  \item a vertical clean equal-ratio curve $C$ has $B_{\mot}(C)\ne0$
  and thus, by (1), carries a Malgrange eigencycle,
  supported arbitrarily close to that curve, in
  $H_{n-1}(M_0;\C)_\lambda$;
  \item a horizontal equal-ratio curve maps onto an algebraic curve
  $\Gamma\subset f^{-1}(0)$ through $0$, and every small non-zero $x\in\Gamma$ has $$
       H_{n-2}(M_x;\C)_\lambda\ne0.       $$
  \end{enumerate}
  \item Every principal Rodrigues--Veys pole of exact order $n-1$ is a root of the Bernstein--Sato polynomial of $f$ at 0.
\end{enumerate}
In particular, the four assertions remain valid for poles of the Hodge
and topological specializations, with ordinary pole order in those settings.
\end{maintheor}
Assertion (1) follows from Theorem~\ref{thm:res-individual-RV-curve} and Corollary~\ref{cor:poles-individual-motivic}, (2) from Corollary~\ref{cor:poles-maximal-motivic} and Theorem~\ref{thm:maximal-order-tube}, (3) from Theorem~\ref{thm:poles-global-RV} and Proposition~\ref{prop:poles-horizontal}, and (4) from Theorem~\ref{thm:bs-high-order}.

\subsection{The structure of the paper} Section~\ref{sec:covers} introduces the preliminaries, and Section~\ref{sec:tubes} the Malgrange eigencycles. Section~\ref{sec:residues} uses non-exactness of residue forms to construct a non-trivial Malgrange eigencycle, hence proving the second part of Theorem~\ref{mainth}. Section~\ref{sec:corner} proves Theorem~\ref{thmmp}, Section~\ref{sec:poles} extends Theorem~\ref{mainth} to motivic zeta-functions, and Section~\ref{sec:bernstein-sato} to Bernstein--Sato polynomials.

\section{Resolution strata, residues, and the topological zeta function}
\label{sec:covers}

We recall the finite covers of resolution strata, define the residue forms
on them, and derive the contribution of an equal-ratio stratum to the
topological zeta function. 

\subsection{Embedded resolution and numerical data}

We work with the polynomial and the resolution of its singularity as in the introduction: $U\subset\A^n_{\C}$ is a Zariski neighborhood of the origin and the resolution
$$
             \rho:X\to U
$$
is projective, is an isomorphism over $U\setminus\{f=0\}$, and has
$\rho^{-1}(0)$ equal to a union of total-transform components. Thus $X$ is smooth and the
reduced divisor
$$
             F^{-1}(0)_{\mathrm{red}}
               =\bigcup_{i\in J}E_i,\quad F:=f\circ\rho,
$$
has simple normal crossings.  Write
\begin{equation}\label{eq:cover-numerical-data}
 \Div(F)=\sum_{i\in J}N_iE_i,
 \qquad
 \Div(\rho^*\omega)=\sum_{i\in J}(\nu_i-1)E_i,
 \qquad
 \omega=dx_1\wedge\cdots\wedge dx_n.
\end{equation}
Components that are strict transforms of components of $\{f=0\}$ are included in
$J$.  If the irreducible factor $f_a$ occurs in $f$ with multiplicity $m_a$, then
its strict transform $E_i$ has $(N_i,\nu_i)=(m_a,1)$.  In particular, the ratio of
a strict component is $1/m_a$.

For a nonempty set $I\subset J$, recall the notation
\begin{equation}\label{eq:Sigma-I}
 E_I=\bigcap_{i\in I}E_i,
 \qquad
 E^\circ_I=E_I\setminus\bigcup_{j\notin I}E_j.
\end{equation}
Thus the $E^\circ_I$ are the open strata of the SNC divisor; in particular,
$E_i^\circ:=E^\circ_{\{i\}}$.  Every
nonempty $E^\circ_I$ is smooth of dimension $n-|I|$.

Near a point of $E^\circ_I$, where $r=|I|$, choose coordinates
$$
     (y_1,\ldots,y_r,z_1,\ldots,z_{n-r})
$$
such that $E_i=(y_i=0)$ for $i\in I$.  
We use $(y,z)$ for these coordinates on $X$ and $(v,z)$ for the
induced coordinates on $N_{E_I^\circ/X}$, where $v_i$ is the linear
function induced by $dy_i$ on each normal fiber.  The $z$-coordinates
are the base coordinates in both systems.
After renumbering the components in
$I$, the two expressions in \eqref{eq:cover-numerical-data} take the form
\begin{align}
 F&=u(y,z)y_1^{N_1}\cdots y_r^{N_r},
       \label{eq:cover-local-F}\\
 \rho^*\omega
   &=b(y,z)y_1^{\nu_1-1}\cdots y_r^{\nu_r-1}
       dy_1\wedge\cdots\wedge dy_r\wedge
       dz_1\wedge\cdots\wedge dz_{n-r},
       \label{eq:cover-local-volume}
\end{align}
where $u$ and $b$ are local units.  These are the coordinates used below for
residue and period calculations.

\subsection{The finite cover of a stratum}

Set
\begin{equation}\label{eq:m-I}
                       m_I=\gcd(N_i\mid i\in I).
\end{equation}
On an open set $V\subset E^\circ_I$ over which coordinates as above have
been chosen, define
\begin{equation}\label{eq:DL-cover-coordinate}
 \{(z,w)\in V\times\CC\mid u(0,z)w^{m_I}=1\}.
\end{equation}
The action $w\mapsto\zeta w$, where $\zeta\in\mu_{m_I}$ and
$\mu_{m_I}$ is the group of $m_I$-th roots of unity, makes this a finite
cyclic torsor over $V$. 
\begin{defin}
These local covers glue into the
\emph{Denef--Loeser cover}:
\begin{equation}\label{eq:DL-cover}
       \widehat E^\circ_I\to E^\circ_I.
\end{equation}
\end{defin}
It is a $\mu_{m_I}$-torsor, note that its total space may be connected.

The same cover has a direct topological interpretation. A normal retraction of a small tubular neighborhood of $E^\circ_{I}$, restricted to the Milnor fiber $F=t$ for a fixed smaller
$t\ne0$, is a locally trivial fibration over $E^\circ_{I}$ (outside a small neighborhood of the infinity). In the normal coordinates $(v,z)$, it becomes the coordinate projection of the set
\begin{equation}\label{eq:leading-coordinate-fiber}
                 u(0,z)\prod_{i\in I}v_i^{N_i}=t.
\end{equation}
Contracting each connected component of every normal fiber over
$z\in E^\circ_I$ to a point gives a cover isomorphic to
\eqref{eq:DL-cover}.

\vspace{1ex}

Suppose now that all components $E_i$ indexed by $I$ have the same ratio
\begin{equation}\label{eq:equal-ratio-on-I}
                       \frac{\nu_i}{N_i}=q\quad(i\in I).
\end{equation}
Write $q=a/d$ in lowest terms.  Since $qN_i=\nu_i\in\Z$, one has
$d\mid N_i$ for every $i\in I$, hence $d\mid m_I$.  Let
\begin{equation}\label{eq:character-q}
 \chi_q:\mu_{m_I}\to\CC,
 \qquad
 \chi_q(e^{2\pi i/m_I})=e^{2\pi iq}.
\end{equation}
This character has order $d$.  
\begin{defin}
The quotient
\begin{equation}\label{eq:effective-cover-general}
 \widehat E^\circ_{I,q}:=\widehat E^\circ_I/\ker\chi_q
      \to E^\circ_I
\end{equation}
is the \emph{$q$-effective cover}. 
\end{defin}
It is a $\mu_d$-torsor,
its total space may be disconnected.
It is important not to replace it by one arbitrarily chosen connected
component: the stabilizer of such a component can have order strictly
smaller than $d$.

In the coordinates of~\eqref{eq:DL-cover-coordinate}, the power
\begin{equation}\label{eq:uq-on-cover}
                 u(0,z)^{-q}=w^{a m_I/d}
\end{equation}
is single-valued on $\widehat E^\circ_{I,q}$: it is invariant under
$\ker\chi_q$ and is multiplied by $\chi_q$ under deck transformations:
$$
 \xi^*\bigl(w^{a m_I/d}\bigr)
       =\chi_q(\xi)w^{a m_I/d},\qquad \xi\in\mu_{m_I}.
$$

\begin{exa}
For a vertical component in the $n=2$ case, the $q$-effective cover is the quotient of
\eqref{eq:DL-cover-coordinate}, with $m_I=N_i$, detected by the character of order $d$.
For a horizontal strict transform of multiplicity $m$, the data are
$(N_i,\nu_i)=(m,1)$, so its ratio is $1/m$ and the cover enumerates the $m$ local roots of $f=t$.
\end{exa}

\subsection{The residue}
\label{subsec:residue-form}

Let $I\subset J$ be nonempty, put $r=|I|$ and $k=n-r$, and assume that
the components containing $E^\circ_I$ have a common ratio
\begin{equation}\label{eq:res-equal-ratio}
                    \frac{\nu_i}{N_i}=q=\frac ad>0
                    \qquad(i\in I),
                    \qquad \gcd(a,d)=1.
\end{equation}
Recall that $m_I=\gcd(N_i:i\in I)$ is divisible by $d$.  The full
finite cover $\widehat E^\circ_I\to E^\circ_I$, the $q$-effective cover
$\widehat E^\circ_{I,q}\to E^\circ_I$, and the character $\chi_q$ were defined in
\eqref{eq:DL-cover}, \eqref{eq:effective-cover-general}, and
\eqref{eq:character-q}.

In this setting, we define the $\chi_q$-equivariant residue form $W_{I,q}$ on $\widehat E^\circ_{I,q}$ as follows. 

In the adapted coordinates of \eqref{eq:cover-local-F} and
\eqref{eq:cover-local-volume}, one has
\begin{align}
 F&=u(y,z)\prod_{i\in I}y_i^{N_i},                                  \notag\\
 \rho^*\omega
  &=b(y,z)\prod_{i\in I}y_i^{qN_i-1}
       \bigwedge_{i\in I}dy_i\wedge dz_1\wedge\cdots\wedge dz_k,
                                                        \label{eq:res-local-pair}
\end{align}
where $u$ and $b$ are local units.
The equalities $\nu_i=qN_i$ give an important cancellation:
\begin{equation}\label{eq:res-log-form}
 F^{-q}\rho^*\omega
 =b(y,z)u(y,z)^{-q}
   \bigwedge_{i\in I}\frac{dy_i}{y_i}
     \wedge dz_1\wedge\cdots\wedge dz_k.
\end{equation}
By \eqref{eq:uq-on-cover}, the power $u(0,z)^{-q}$ is single-valued on the cover $\widehat E^\circ_{I}$, where it equals $w^{am_I/d}$. So the iterated
logarithmic
residue of (\ref{eq:res-log-form}) is a single-valued form
\begin{equation}\label{eq:res-local-residue-general}
 \widehat W_{I,q}
 =b(0,z)w^{a m_I/d}\,dz_1\wedge\cdots\wedge dz_k
\end{equation}
on this cover. It is invariant under $\ker\chi_q$, and hence descends to a
nowhere-vanishing holomorphic $k$-form $W_{I,q}$ on $\widehat E^\circ_{I,q}$.
With the convention \eqref{eq:character-q}, we have
\begin{equation}\label{eq:res-equivariance}
                 \zeta^*\widehat W_{I,q}
                    =\chi_q(\zeta)\widehat W_{I,q}.
\end{equation}

The residue is independent of the adapted coordinates. Thus the local expressions \eqref{eq:res-local-residue-general}
define $W_{I,q}$ globally; an ordering of $I$ fixes only its overall sign.

The residue form is holomorphic and nowhere vanishing on $\widehat E^\circ_{I,q}$, but may have zeroes and poles at the boundary of a compactification of its domain. 
We characterize the orders of these zeroes and poles as follows.

Suppose that the closure of $E^\circ_I$ meets a further component $E_j$
transversely along a divisor $B$, and put
\begin{equation}\label{eq:res-general-defect}
                            \alpha_j=\nu_j-qN_j.
\end{equation}
In a coordinate $z_1$ transverse to $B$ inside the closure of the stratum,
the residue is a unit times $z_1^{\alpha_j-1}dz_1$ wedged with a nonvanishing holomorphic 
$(k-1)$-form.  If a component of the $q$-effective cover has ramification
index $e$ along $B$, substitution $z_1=\widetilde z_1^e$ gives
\begin{equation}\label{eq:res-general-order}
                    \ord_{\widetilde B}W_{I,q}
                                  =e\alpha_j-1.
\end{equation}

\subsection{The topological zeta function and strata contributions at 0}
\label{subsec:topological-contributions}

We use local zeta functions at the origin; at another point $x\in f^{-1}(0)$,
the same definitions apply after translating $x$ to $0$.
The resolution formula for the local topological zeta function is
\begin{equation}\label{eq:topological-resolution-formula}
 Z_{\mathrm{top}}(f;s)
  =\sum_{\varnothing\ne I\subset J}
       \chi_c(E^\circ_I\cap\rho^{-1}(0))
       \prod_{i\in I}\frac1{\nu_i+sN_i},
\end{equation}
where $\chi_c$ is the compactly supported Euler characteristic
\cite{DenefLoeser1992}. This rational function is independent of the
resolution. Its poles belong to the set $\{-\nu_i/N_i\}$, but cancellation
between the summands may remove a candidate pole or lower its order.

Fix $q>0$, a nonempty $I\subset J$, and a connected component $Y$ of $E_I$,
where $r=|I|$ and
$\nu_i/N_i=q$ for every $i\in I$. Thus $\dim Y=n-r$. Call $Y$ \emph{clean}
if no further ratio-$q$ component meets it. Put $\alpha_j=\nu_j-qN_j$.
The intersections $Y\cap E_{I\cup K}^\circ$, for $K\subset J\setminus I$,
partition $Y$.
Grouping the terms of \eqref{eq:topological-resolution-formula} lying on
$Y$ gives
\begin{equation}\label{eq:topological-stratum-package}
 \frac1{(s+q)^r\prod_{i\in I}N_i}
 \sum_{K\subset J\setminus I}
       \chi_c(Y\cap E_{I\cup K}^\circ\cap\rho^{-1}(0))
       \prod_{j\in K}\frac1{\alpha_j+(s+q)N_j}.
\end{equation}
Here and below empty strata are omitted and an empty product is $1$.
For clean $Y$, the sum in this expression is regular at $s=-q$.
\begin{defin}
The \emph{topological contribution} of $Y$ at the origin is its value there:
\begin{equation}\label{eq:general-Btop}
 B_{\mathrm{top}}(Y)
   =\sum_{K\subset J\setminus I}
       \chi_c(Y\cap E_{I\cup K}^\circ\cap\rho^{-1}(0))\prod_{j\in K}\frac1{\alpha_j}.
\end{equation}
\end{defin}
Thus its contribution to the coefficient of $(s+q)^{-r}$ is
$B_{\mathrm{top}}(Y)/\prod_{i\in I}N_i$. This is an individual contribution;
contributions from different components may cancel in the whole zeta
function. If $Y$ is not clean, \eqref{eq:topological-stratum-package}
can have higher-order poles, and \eqref{eq:general-Btop} does not apply.

\begin{exa}
For a complete clean curve $C\subset\rho^{-1}(0)$, let $P=C\setminus C^\circ$ and let
$\alpha_p$ be the defect of the component transverse to $C$ at $p$.
\begin{equation}\label{eq:res-Btop}
 B_{\mathrm{top}}(C)
   =\chi_c(C^\circ)+\sum_{p\in P}\frac1{\alpha_p}
   =2-2g(C)-|P|+\sum_{p\in P}\frac1{\alpha_p}.
\end{equation}
For a point $Y=\{p\}\subset\rho^{-1}(0)$ where $n$ ratio-$q$ components
meet, there are no further components and $B_{\mathrm{top}}(Y)=1$.
\end{exa}

The next two sections explain how non-zero zeta contributions of curve strata over 0 produce respective eigencycles in the homology of the  Milnor fiber at 0. How the strata outside 0 produce eigencycles at their nearby points is the easy part; we recall the respective standard fact here.
\begin{utver}[Nearby eigencycles]
\label{prop:poles-horizontal}
Let $n\ge2$, let $I\subset J$ have $|I|=n-1$, and suppose that
$\nu_i/N_i=q$ for every $i\in I$.  Let $C$ be a horizontal
irreducible component of $E_I$ which meets $\rho^{-1}(0)$, and let
$\Gamma=\rho(C)$ be its image, which is closed because $\rho$ is proper. Put
$\lambda=e^{-2\pi iq}$. Then, for all sufficiently small $x\in\Gamma\setminus\{0\}$, we have \begin{equation}\label{eq:poles-horizontal-eigenspace}
                   H_{n-2}(M_x;\C)_\lambda\ne0.
\end{equation}
For every such $x$ and every lift
$y\in C\cap\rho^{-1}(x)$, the
class in \eqref{eq:poles-horizontal-eigenspace} may be supported
arbitrarily close to $y$ in the inverse image of the Milnor fiber.  Moreover,
\begin{equation}\label{eq:poles-horizontal-q}
                              q=\frac1d
                  \qquad\text{for some }d\in\Z_{>0}.
\end{equation}
\end{utver}

\section{Leray tubes and Malgrange eigencycles}
\label{sec:tubes}

A compact singular cycle $\gamma$ in the cover $\widehat E_I^\circ$ gives
a compact cycle in the nearby fiber $F=t$: its \emph{Leray tube}.
The Leray tube:

-- is supported in an arbitrarily small neighborhood of the projection of the support $|\gamma|\subset \widehat E_I^\circ$, for sufficiently small $t$;

-- is transformed by the monodromy of $F$ when $\gamma$ is transformed by the deck
transformation of the cover;

-- for an equal-ratio stratum $E_I^\circ$, has the Gelfand--Leray period whose leading
asymptotics in $t$ is given by the residue period of $\gamma$.

In particular, for a prescribed character $\chi$ of the deck transformation group $G$, the $\chi$-equivariant symmetrization $\sum_{g\in G} \chi(g)^{-1} \cdot g_* \gamma$ has a Leray tube whose homology class is a monodromy
eigenvector (or zero).  This tube is called the \emph{Malgrange eigencycle}
of $\gamma$.  

We introduce these notions and prove their properties. Our motivation is to prove the monodromy conjecture by constructing relevant Malgrange eigencycles, whose non-vanishing is certified by their Gelfand-Leray periods, and thus yields the sought monodromy eigenvalue.

\subsection{Leray tubes}
\label{subsec:normal-tori}

Put $r=|I|$.  The normal bundle of $E_I^\circ$ splits into the normal
lines of the divisors $E_i$, $i\in I$.  At $z\in E_I^\circ$, the
coordinates $(y,z)$ of~\eqref{eq:cover-local-F} induce normal coordinates
$(v,z)$ by $v_i=dy_i$.
In these coordinates, the leading normal term of $F=f\circ\rho$ is
the function
\begin{equation}\label{eq:normal-leading-tube}
 F_0(z;v_1,\ldots,v_r)=u(0,z)\prod_{i\in I}v_i^{N_i}.
\end{equation}
For $t\ne0$, the fibers of the normal level $F_0=t$ over $E_I^\circ$
are cosets of
$$
           \{(\lambda_i)\in\CC^I:\ \prod_i\lambda_i^{N_i}=1\}.
$$
This subgroup has $m_I$ connected components, each isomorphic to
$\CC^{r-1}$.

Fix smooth Hermitian norms on the normal line bundles and put
\begin{equation}\label{eq:equal-norm-tori}
 \mathcal T_t=\{v\in N_{E_I^\circ/X}:\ F_0(v)=t,\quad
                         \|v_1\|=\cdots=\|v_r\|\}.
\end{equation}
Its projection to $E_I^\circ$ is a smooth locally trivial fibration,
whose fiber has $m_I$ connected components, each a compact real
$(r-1)$-torus onto which the corresponding complex torus component of $F_0=t$ 
retracts.  

As recalled in Section~\ref{sec:covers}, contracting the connected
components of the normal fibers gives $\widehat E_I^\circ$. The same component cover for
$\mathcal T_t$ gives the bundle
\begin{equation}\label{eq:leading-mapping-torus}
 p_t:\mathcal T_t\longrightarrow\widehat E_I^\circ
 \qquad\hbox{with fiber }(S^1)^{r-1}.
\end{equation}
Its fibres are oriented by an orientation of the lattice
\begin{equation}\label{eq:kernel-lattice-general}
             \{(c_i)\in\Z^I:\ \sum_iN_ic_i=0\}.
\end{equation}

Fix a smooth diffeomorphism $\Phi$ from a neighborhood of the zero
section in the normal bundle to a neighborhood of $E_I^\circ$ in $X$,
fixing $E_I^\circ$ and inducing the identity on its complex normal bundle,
such that
\begin{equation}\label{eq:exact-normal-identification}
                              F\circ\Phi=F_0.
\end{equation}
It identifies the normal level $F_0=t$, restricted to its domain, with
the corresponding part of $F=t$.  The compact tori $\mathcal T_t$
therefore are sent to compact tori in the actual Milnor fiber.

\begin{defin}
The \emph{Leray tube} $\Tub_t(\gamma)$ of a finite singular chain
$\gamma$ on $\widehat E_I^\circ$ is the image under $\Phi$ of its
oriented inverse image under $p_t$.
\end{defin}
Concretely, pull the torus bundle back to each simplex of $\gamma$ and
triangulate compatibly along faces.  Orient the base before the fiber.
This gives a chain of dimension $\dim\gamma+r-1$ and
$$
                    \partial\Tub_t(\gamma)=\Tub_t(\partial\gamma).
$$
Thus cycles give compact closed cycles.    Changing the
auxiliary choices through a family gives homologous tubes.
\begin{exa}
A point gives one normal torus;
a smooth loop gives a torus bundle over that loop.
\end{exa}

\subsection{Monodromy and periods of Leray tubes}
We now describe an explicit monodromy representative of $F_0$, which through $\Phi$ induces the monodromy of $F$.
Write $v=(v_1,\ldots,v_r)$.  Along $t(s)=te^{2\pi is}$ the rotations
\begin{equation}\label{eq:monomial-transport}
 R_s(z,v)=(z,e^{2\pi is/\sum_{j\in I}N_j}v),
             \qquad 0\le s\le1,
\end{equation}
carry $F_0=t$ to $F_0=t(s)$ and preserve the norms, hence carry
$\mathcal T_t$ to $\mathcal T_{t(s)}$.  Their endpoint preserves $\mathcal T_t$ and
induces the deck generator $\zeta_I=e^{2\pi i/m_I}$ on
$\widehat E_I^\circ$.  

These are the rotations in A'Campo's construction at the simplex
coordinates $N_i/\sum_{j\in I}N_j$; his $i$-th angular speed is the
$i$-th simplex coordinate divided by $N_i$.  See~\cite[\S2(1),
p.~240]{ACampo1975} for the construction and~\cite[\S2(2)--(3),
pp.~241--242, Theorem~5]{ACampo1975} for its comparison with Milnor
monodromy.  Through $\Phi$ these rotations give geometric transport
near our tori; for their homology classes, the comparison also follows
directly from the continuous family of cycles in the proof below.

For $x\in U\cap\{f=0\}$, fix a sufficiently small closed Milnor ball
$B_\epsilon(x)\subset U$.  We retain the local notation
$$
 M_{x,t}=B_\epsilon(x)\cap f^{-1}(t)
 \simeq \rho^{-1}(B_\epsilon(x))\cap F^{-1}(t).
$$

\begin{lemma}\label{lem:tube-support-monodromy}
Let $\gamma$ be a finite singular cycle on $\widehat E_I^\circ$.

(1) For every neighborhood $V$ of its projected support in $X$, the
tube $\Tub_t(\gamma)$ is supported in $V\cap\{F=t\}$ for all
sufficiently small $t$.  If the projected support of $\gamma$ lies in
$\rho^{-1}(x)$, the tube defines a compact closed cycle in the Milnor fiber $M_{x,t}$.

(2) In this case,
\begin{equation}\label{eq:monodromy-tube-relation}
 \M_x[\Tub_t(\gamma)]=[\Tub_t((\zeta_I)_*\gamma)].
\end{equation}
Along $t>0$, the tube classes are transported by the Milnor fibration.
\end{lemma}
\begin{proof}
On a sector of the punctured complex $t$-plane, fix one branch of
$\log t$ and use it to define $t^{1/\sum_iN_i}$ and $t^{1/m_I}$.
Then
$$
 p_t(z,t^{1/\sum_iN_i}v)=p_1(z,v),\qquad
 \Phi(z,t^{1/\sum_iN_i}v)\longrightarrow z\quad(t\to0),
$$
uniformly on the compact set $p_1^{-1}(|\gamma|)$.

Suppose that the projected support lies in $\rho^{-1}(x)$.
The rotations $R_s$ in~\eqref{eq:monomial-transport} preserve the
normal norms defining $\mathcal T_t$ in~\eqref{eq:equal-norm-tori},
so for small $|t|$ the chains
$$
 C_s=\rho\Phi R_s(p_t^{-1}\gamma),\qquad 0\le s\le1,
$$
form a continuous family of cycles in $M_{x,te^{2\pi is}}$.
Choose a trivialization $P_s:M_{x,t}\longrightarrow M_{x,te^{2\pi is}}$
of the Milnor fibration pulled back to $[0,1]$, with $P_0=\mathrm{id}$.
The cycles $P_s^{-1}C_s$ are homotopic, hence
$$
 [C_1]=(P_1)_*[C_0]=\M_x[C_0].
$$
Since $p_tR_1=\zeta_Ip_t$ and $R_1$ preserves the oriented torus fibers,
$[C_1]=[\Tub_t((\zeta_I)_*\gamma)]$ under the identification with the
Milnor fiber.  This gives~\eqref{eq:monodromy-tube-relation}.
The same argument applied to the scalings
$(z,v)\mapsto(z,t^{1/\sum_iN_i}v)$ along $t>0$ proves the last assertion.
\end{proof}

The Gelfand--Leray form $\rho^*\omega/dF$ on the fiber $F^{-1}(t)$ is the
restriction of any local $(n-1)$-form $\eta$ satisfying
$$
                         dF\wedge\eta=\rho^*\omega.
$$
Two such choices have the same restriction to the fiber, so the form is
intrinsic and closed.

\begin{lemma}[Residue periods and tube periods]\label{lem:tube-period}
Assume $\nu_i=qN_i$ for every $i\in I$, and put $k=n-r$.
For every finite piecewise smooth singular $k$-chain $\gamma$ on the full
cover $\widehat E_I^\circ$, with complex coefficients,
\begin{equation}\label{eq:general-tube-period}
 \int_{\Tub_t(\gamma)}\frac{\rho^*\omega}{dF}
   =\pm\frac{(2\pi i)^{r-1}}{m_I}t^{q-1}
       \left(\int_\gamma\widehat W_{I,q}+o(1)\right),
\end{equation}
as functions of $t\to 0$. Here $\widehat W_{I,q}$ is the residue from
Subsection~\ref{subsec:residue-form}, and the sign depends only on the
orientation conventions (irrelevant for this work).  More generally, for every holomorphic
function $h$ near the image of the support downstairs,
\begin{equation}\label{eq:tube-all-numerators}
 \int_{\Tub_t(\gamma)}\frac{\rho^*(h\omega)}{dF}
 =\pm\frac{(2\pi i)^{r-1}}{m_I}t^{q-1}
   \left(\int_\gamma(h\circ\rho)\widehat W_{I,q}+o_h(1)\right).
\end{equation}
The function $h\circ\rho$ in the right-hand integral is restricted to
the stratum.  Piecewise smooth means that, after a finite subdivision,
each simplex parametrization is smooth up to its faces.  Both asymptotics
are uniform for the complex fiber value $t$ tending to $0$ within
any fixed angular sector
$$
 \{t\in\C:0<|t|<\delta,\quad\alpha\le\arg t\le\beta\},
 \qquad \beta-\alpha<2\pi.
$$
All powers of $t$ use the same branch of $\log t$ on this sector
of the punctured complex $t$-plane. \end{lemma}

If the projected support lies over $x$, the restriction of $h\circ\rho$
is the constant $h(x)$.  We shall use this case of
\eqref{eq:tube-all-numerators} in Section~\ref{sec:bernstein-sato}.

\begin{proof}
Locally, up to the orientation sign,
\begin{equation}\label{eq:general-GL-form}
 \frac{\rho^*\omega}{dF}=t^{q-1}
 \frac{b(y,z)u(y,z)^{-q}}
      {N_r+y_r\partial_{y_r}\log u(y,z)}
 \bigwedge_{i=1}^{r-1}\frac{dy_i}{y_i}\wedge
 dz_1\wedge\cdots\wedge dz_k.
\end{equation}
By subdivision and linearity, consider a smooth simplex
$\sigma:\Delta^k\to\widehat E_I^\circ$, $\sigma(s)=(z(s),w(s))$,
lying over an adapted chart.  Choose a smooth section
$a(s)=(a_1(s),\ldots,a_r(s))$ of $\sigma^*\mathcal T_1$; thus
$\prod_i a_i(s)^{N_i/m_I}=w(s)$.  Put
$$
 K=\{(\tau_1,\ldots,\tau_r)\in(S^1)^r:
                         \prod_i\tau_i^{N_i/m_I}=1\}.
$$
Since the integers $N_i/m_I$ are relatively prime, $K$ is a connected
$(r-1)$-torus.  The tube is parametrized by
$$
 \begin{gathered}
 \Psi_t:\Delta^k\times K\longrightarrow F^{-1}(t)\subset X,\\
 \Psi_t(s,\tau)=\Phi\bigl(z(s),
 t^{1/\sum_jN_j}(a_1(s)\tau_1,\ldots,a_r(s)\tau_r)\bigr).
 \end{gathered}
$$
In the same adapted chart $(y,z)$ of $X$, write
$(Y_t,Z_t)=(y,z)\circ\Psi_t$.
The assumptions on $\Phi$ give
$$
 Y_{t,i}=t^{1/\sum_jN_j}a_i(s)\tau_i
                     +O_{C^1}(|t|^{2/\sum_jN_j}),\qquad
 Z_t=z(s)+O_{C^1}(|t|^{1/\sum_jN_j}),
$$
where the estimates include first derivatives on $\Delta^k\times K$.
Using the branch
$$
 u(Y_t,Z_t)^{-q}=t^{-q}\prod_iY_{t,i}^{\nu_i}
                 \longrightarrow w(s)^{qm_I}
$$
in~\eqref{eq:general-GL-form}, and~\eqref{eq:res-local-residue-general},
we obtain uniformly
$$
 t^{1-q}\Psi_t^*\frac{\rho^*\omega}{dF}
 \longrightarrow
 \pm\frac1{N_r}\,\sigma^*\widehat W_{I,q}\wedge
                    \bigwedge_{i=1}^{r-1}\frac{d\tau_i}{\tau_i}.
$$
Here the terms with $k+1$ base differentials vanish.  The projection
$K\to(S^1)^{r-1}$ onto the first $r-1$ coordinates is a covering of
degree $N_r/m_I$, so
$$
 \int_K\bigwedge_{i=1}^{r-1}\frac{d\tau_i}{\tau_i}
       =\pm(2\pi i)^{r-1}\frac{N_r}{m_I}.
$$
Integration proves~\eqref{eq:general-tube-period}.  For $r=1$, $K$
is a point and the empty wedge and its integral are $1$.
Multiplying by $h\circ\rho\circ\Psi_t\to(h\circ\rho)(0,z(s))$
proves~\eqref{eq:tube-all-numerators}.
\end{proof}

\subsection{Malgrange eigencycles and the period certificate}

Assume that the components containing $E_I^\circ$ have common ratio
$q=a/d$ in lowest terms.  Recall that the covers
$$
 \widehat E_I^\circ\longrightarrow E_I^\circ,\qquad
 \widehat E_{I,q}^\circ\longrightarrow E_I^\circ
$$
have degrees $m_I$ and $d$, respectively.  The quotient map
$$
 p:\widehat E_I^\circ\longrightarrow\widehat E_{I,q}^\circ
       =\widehat E_I^\circ/\ker\chi_q
 \qquad\hbox{has degree }|\ker\chi_q|=m_I/d.
$$
The monodromy of $F_0$ (given by the endpoint $R_1$ of~\eqref{eq:monomial-transport}) induces the deck
generator $\zeta_I$ upstairs and its image $\bar\zeta_I$ of order $d$
downstairs:
$$
 p_tR_1=\zeta_Ip_t,\qquad
 p\zeta_I=\bar\zeta_Ip. $$
The character $\chi_q$ factors through this quotient and sends
$\bar\zeta_I$ to $e^{2\pi ia/d}$.

Let $\gamma$ be a finite singular chain on the $q$-effective cover.
Choose one lift of each simplex to the full cover, retaining its
coefficient, and denote the resulting chain by $\widetilde\gamma$.
Define its normalized $\chi_q$-equivariant symmetrization by
\begin{equation}\label{eq:normalized-transfer}
 \Gamma=\frac1{m_I}\sum_{\xi\in\mu_{m_I}}
                  \chi_q(\xi)^{-1}\xi_*\widetilde\gamma.
\end{equation}
This chain is independent of the choices, and satisfies  $$
 \partial\Gamma(\gamma)=\Gamma(\partial\gamma),\qquad
 (\zeta_I)_*\Gamma=e^{2\pi iq}\Gamma.
$$
In particular, $\Gamma$ is a cycle if $\gamma$ is.

\begin{rem}
If $p^!\gamma$ is the sum of all lifts of every simplex of $\gamma$,
with their coefficients, then
$$
 \Gamma=\frac1{m_I|\ker\chi_q|}
       \sum_{\xi\in\mu_{m_I}}\chi_q(\xi)^{-1}\xi_*p^!\gamma.
$$
\end{rem}

If $\gamma$ is a cycle whose support projects to $\rho^{-1}(x)$,
where $x\in U\cap\{f=0\}$, Lemma~\ref{lem:tube-support-monodromy} gives
\begin{equation}\label{eq:tube-eigen-relation}
              \M_x[\Tub_t(\Gamma)]
                    =e^{2\pi iq}[\Tub_t(\Gamma)].
\end{equation}

\begin{defin}
For a finite singular cycle $\gamma$, its \emph{Malgrange eigencycle}
for the character $\chi_q$ is the Leray tube of the cycle $\Gamma$
defined by~\eqref{eq:normalized-transfer}.
\end{defin}
Note that $\gamma\ne 0$ may have Malgrange eigencycle vanishing in the Milnor fiber homology. We now learn to certify it does not happen, if an appropriate residue period of $\gamma$ is not 0.

Put $k=n-r$.  Suppose that $\gamma$ is a finite piecewise smooth $k$-cycle with
\begin{equation}\label{eq:nonzero-residue-period}
                         \int_\gamma W_{I,q}\ne0.
\end{equation}
Equations~\eqref{eq:normalized-transfer} and~\eqref{eq:res-equivariance}
give
\begin{equation}\label{eq:Fourier-period-preserved}
                 \int_\Gamma\widehat W_{I,q}
                    =\int_\gamma W_{I,q}\ne0.
\end{equation}
The period lemma now gives the following criterion.

\begin{utver}[Malgrange eigencycle criterion]
\label{prop:general-Fourier-Leray}
Assume that the components containing $E^\circ_I$ have common ratio
$q>0$.  Let $W_{I,q}$ be the residue form on the $q$-effective cover,
characterized by $p^*W_{I,q}=\widehat W_{I,q}$ and the local
formula~\eqref{eq:res-local-residue-general}.  If it has a finite
piecewise smooth $k$-cycle $\gamma$ with nonzero period, where
$k=n-|I|$, then
$$
       e^{-2\pi iq}
$$
is an eigenvalue of geometric monodromy on
$H_{n-1}(M_x;\C)$ whenever the image of $|\gamma|$ under
$\widehat E_{I,q}^\circ\to E_I^\circ\subset X$ lies in
$\rho^{-1}(x)$, for $x\in U\cap\{f=0\}$.
For every neighborhood $V$ of this image in $X$, the eigenspace contains
a nonzero class (particularly, the coefficientwise complex conjugate to the Malgrange eigencycle of $\gamma$) in the image of
$$
 H_{n-1}\bigl(V\cap\rho^{-1}(M_{x,t});\C\bigr)
       \to H_{n-1}(M_x;\C)
$$
for all sufficiently small $t\ne0$.
\end{utver}

\begin{proof}
Lemma~\ref{lem:tube-period} and~\eqref{eq:Fourier-period-preserved}
give a nonzero Gelfand--Leray period of $\Tub_t(\Gamma)$ for every
sufficiently small $t$.  This remains true when $q>1$, although the
period then tends to zero.  Since the form is closed, a homologically
zero tube would bound a finite piecewise smooth singular chain and Stokes'
theorem would make the period zero.
Thus the tube represents a nonzero ordinary homology class.
Its eigenvalue and support follow from
\eqref{eq:tube-eigen-relation} and
Lemma~\ref{lem:tube-support-monodromy}; coefficientwise complex
conjugation gives $e^{-2\pi iq}$ with the same support.
\end{proof}
\begin{exa}
In case $n=2$, if $E_I$ is a vertical component, formula~\eqref{eq:general-tube-period}
has $r=1$ and becomes
\begin{equation}\label{eq:plane-curve-tube-period}
 \int_{\Tub_t(\Gamma)}\frac{\rho^*\omega}{dF}
  =\pm\frac{1}{N_i}t^{q-1}
      \left(\int_\Gamma\widehat W_{i,q}+o(1)\right).
\end{equation}

There is no factor $2\pi i$, because there is no continuous normal
circle.  At a plane full corner, the base is a point and the formula
has one normal factor $2\pi i$.  \end{exa}
\begin{rem}\label{rem:degree-of-tube}
The homology dimension $n-1$ comes from
$$
        \dim_{\mathbf R}(\text{base cycle})
        +\dim_{\mathbf R}(\text{normal torus})
           =k+(|I|-1)=n-1.
$$
For a horizontal stratum one instead takes a normal slice; its ambient
complex dimension is smaller, and the same formula gives the corresponding
lower homological degree.  That distinction will be made explicit when
horizontal equal-ratio curves are treated.

When $n=2$, a normal slice to a horizontal curve $E_I$ is one-dimensional and its Milnor
fiber is a finite set.  The resulting Malgrange eigencycle lies in $H_0$ homology;
it lies in reduced $\widetilde H_0$ precisely for a nontrivial character.
\end{rem}

\subsection{Illustration: poles of maximal order}

We illustrate the use of Malgrange eigencycles by proving the monodromy eigenvalue exists for a
pole of the largest possible order.

\begin{theor}[Maximal-order poles]
\label{thm:maximal-order-tube}
Let $x\in f^{-1}(0)$ and $q>0$. Suppose that $s=-q$ is a pole of order
$n$ of the local topological zeta function (or, more generally, a principal Hodge
pole). Then $e^{-2\pi iq}$ is an eigenvalue of $\M_x$ on
$H_{n-1}(M_x;\C)$.  More precisely, there is a point
$p\in\rho^{-1}(x)$ at which $n$ ratio-$q$ components meet, and a nonzero
eigenclass can be represented in every sufficiently small neighborhood of
$p$ in the inverse image of the Milnor fiber.
\end{theor}

\begin{proof}
At most $n$ components of an SNC divisor on $X$ can meet.  The coefficient
of $(s+q)^{-n}$ in \eqref{eq:topological-resolution-formula} is
\begin{equation}\label{eq:maximal-top-coefficient}
 \sum_{\substack{p\in\rho^{-1}(x)\\
       p\text{ lies on }n\text{ components }E_i\\
       \nu_i/N_i=q}}
              \frac1{\prod_{E_i\ni p}N_i}.
\end{equation}
Every summand is positive.  Hence a topological pole of order $n$ is
equivalent to the existence of at least one point in this sum.
(For a principal Hodge pole of order $n$, the same
nonempty set of points is given by Corollary~\ref{cor:poles-maximal-motivic}.)  

Take adapted coordinates at $p$:
\begin{align*}
 F&=u(y)y_1^{N_1}\cdots y_n^{N_n},\\
 \rho^*\omega
   &=b(y)y_1^{qN_1-1}\cdots y_n^{qN_n-1}
          dy_1\wedge\cdots\wedge dy_n,
\end{align*}
with $u(0)b(0)\ne0$.  Put $m=\gcd(N_1,\ldots,N_n)$.  If $q=a/d$ is
reduced, then $d\mid m$.  The $q$-effective cover of the point is a
$\mu_d$-torsor.  Let $\Gamma$ be the normalized $\chi_q$-equivariant symmetrization
\eqref{eq:normalized-transfer} of any one of its points.  Formula
\eqref{eq:general-tube-period}, now with $r=n$ and $k=0$, becomes
\begin{equation}\label{eq:maximal-torus-period}
 \int_{\Tub_t(\Gamma)}\frac{\rho^*\omega}{dF}
   =\pm\frac{(2\pi i)^{n-1}}m
       t^{q-1}\bigl(b(0)u(0)^{-q}+o(1)\bigr)\ne0
\end{equation}
for all sufficiently small $t\ne0$.  Here $\Tub_t(\Gamma)$ is a finite,
character-weighted sum of compact real $(n-1)$-tori, hence the Malgrange
eigencycle of the corner.  It is an absolute cycle, and its
period proves that its homology class is nonzero, and
\eqref{eq:monodromy-tube-relation} gives the required character.  \end{proof}

The proof also clarifies why the order-$n$ case is particularly direct:
the stratum is a point, the character residue is the nonzero scalar
$b(0)u(0)^{-q}$, and nonexactness questions do not arise.  For order
$n-1$, the corresponding stratum is a curve; the next section studies the exactness of the 
meromorphic one-form that replaces this scalar.
\section{Nonexactness of curve residues}
\label{sec:residues}

We work in the setting of Section~\ref{sec:covers} with covers and residue forms introduced there.
If an equal-ratio stratum of the exceptional divisor is a vertical curve, we prove nonexactness of its residue form under a numerical hypothesis below. Nonexactness implies there is a loop with a non-zero integral of the residue form, which certifies that the Malgrange eigencycle over this loop is non-zero in the homology of the Milnor fiber, as proved in Section~\ref{sec:tubes}.

\subsection{Vertical equal-ratio curves}
\label{subsec:complete-equal-ratio-curves}

Let $C\subset\rho^{-1}(0)$ be a
component of the smooth curve $E_1\cap\cdots\cap E_{n-1}$, assume
$$
                 \frac{\nu_1}{N_1}=\cdots=
                 \frac{\nu_{n-1}}{N_{n-1}}=q, \qquad  C^\circ=C\setminus P,
$$
where $P$ is the intersection of $C$ with the rest of the components of the hypersurface $F=0$ (including the strict transform).
At $p\in P$ there is one transverse component, with data $(N_p,\nu_p)$; recall that its {\it defect} is
\begin{equation}\label{eq:res-curve-defect}
                            \alpha_p=\nu_p-qN_p.
\end{equation}
Thus $C$ is {\it clean} in the sense of Subsection~\ref{subsec:topological-contributions}
if and only if $\alpha_p\ne0$ for all $p\in P$.

Retain the notation $q=a/d$.  Choose a connected component
$\widetilde C^\circ$ of the $q$-effective cover,
let $M$ be its degree over $C^\circ$, and let $\widetilde C$ be its smooth
compactification.  Thus $M\mid d$.

\begin{lemma}[Ramification of the cover and the divisor of the residue]
\label{lem:res-curve-inertia}
For every $p\in P$, one has $(d/M)\mid N_p$.  The ramification index at each
point $\widetilde p\in\widetilde C$ over $p$ equals
\begin{equation}\label{eq:res-curve-inertia-formula}
 e_p=\frac d{\gcd(d,N_p)}=(\mbox{the denominator of }\alpha_p)\mid M,
\end{equation}
and the residue differential $W=W_{I,q}|_{\widetilde C}$ satisfies
\begin{equation}\label{eq:res-curve-order-formula}
                         \ord_{\widetilde p}W
                                  =e_p\alpha_p-1.
\end{equation}
In particular, $M\alpha_p\in\Z$ for every $p\in P$.
It has no zero or pole over $C^\circ$.
\end{lemma}

\begin{proof}
Let $z$ be a coordinate on $C$ centered at $p$.  By
\eqref{eq:DL-cover-coordinate}, the cover over a punctured disc is given by
$$
 u(z)z^{N_p}w^{m_I}=1,\qquad u(0)\ne0.
$$
Choose a holomorphic $d$-th root of $u(z)$.  The coordinate
$v=u(z)^{1/d}w^{m_I/d}$ on the $q$-effective cover satisfies
$$
 v^d=z^{-N_p}.
$$
Continuation around $z=0$ sends
$$
                 v\longmapsto e^{-2\pi iN_p/d}v.
$$
The order of this permutation, and the length of each of its orbits,
is $e_p=d/\gcd(d,N_p)$.  This is the ramification index of the
compactified cover.  The chosen connected component has $M$ points
in each fiber and is preserved by continuation, so $e_p\mid M$.
Since $q=a/d$ is in lowest terms, $e_p$ is also the denominator of
$\alpha_p=\nu_p-aN_p/d$.

On a sector of that punctured disc, choose a branch of the rational
power of $z$.  Formula~\eqref{eq:res-local-residue-general} expresses
$W$ as the pullback of $h(z)z^{\alpha_p-1}dz$, where $h$ is holomorphic
and $h(0)\ne0$.  Under the covering map
$\widetilde C\to C$, locally $z=\widetilde z^{e_p}$, this gives
$$
 W=e_p h(\widetilde z^{e_p})
           \widetilde z^{e_p\alpha_p-1}d\widetilde z,
$$
which proves~\eqref{eq:res-curve-order-formula}.
Its integral order gives $e_p\alpha_p\in\Z$; since $e_p\mid M$,
we obtain $M\alpha_p\in\Z$.  Thus $d\mid MN_p$, or
$(d/M)\mid N_p$, as claimed.
Formula~\eqref{eq:res-local-residue-general} also proves that $W$
has no zero or pole over $C^\circ$.
\end{proof}

Assume that $C$ is clean. Its topological contribution
$B_{\mathrm{top}}(C)$ is given by \eqref{eq:res-Btop}.
The following is well known (see \cite{Veys1998}).
\begin{lemma}[Equal-ratio adjunction]
\label{lem:res-equal-ratio-adjunction}
One has
\begin{equation}\label{eq:res-adjunction-identities}
 \sum_{p\in P}(\alpha_p-1)=2g(C)-2,
 \qquad
 B_{\mathrm{top}}(C)=
       \sum_{p\in P}\left(\frac1{\alpha_p}-\alpha_p\right).
\end{equation}
\end{lemma}

A puncture with $\alpha_p=1$ is \emph{neutral} in two compatible senses.
First, filling it back into $C^\circ$ increases $\chi_c(C^\circ)$ by $1$
and removes the summand $1/\alpha_p=1$, so it does not change
$B_{\mathrm{top}}(C)$.
Second, 
$$
                         qN_p=\nu_p-1\in\Z
$$
implies $d\mid N_p$. Hence the character cover extends without ramification
across $p$, and \eqref{eq:res-curve-order-formula} gives
$\ord_{\widetilde p}W=0$, so the residue differential extends without a zero
or pole. 

\subsection{The norm of the residue antiderivative}
\label{subsec:res-exactness}
For the compactified connected cover
$$
                         \pi:\widetilde C\to C
$$
of degree $M$, suppose that the 1-form on the covering space is \emph{exact}:
$$W=dG_0$$ 
for a meromorphic function $G_0$ on $\widetilde C$.
Assume that $M>1$ (the trivial case $M=1$ will be settled at the end). The deck group of $\widetilde C\to C$ is the
subgroup $\mu_M\subset\mu_d$ preserving the chosen connected component.
We use $\chi_q$ both for its descended character on $\mu_d$ and its
restriction to $\mu_M$. This restriction is faithful and satisfies
$$
                              g^*W=\chi_q(g)W.
$$
The normalized $\chi_q$-equivariant symmetrization of an antiderivative is
\begin{equation}\label{eq:res-equivariant-primitive}
 G=\frac1M\sum_{g\in\mu_M}\chi_q(g)^{-1}g^*G_0.
\end{equation}
It satisfies
\begin{equation}\label{eq:res-equivariant-primitive-properties}
                  dG=W,\qquad g^*G=\chi_q(g)G.
\end{equation}
\begin{rem}

This normalization is independent of the additive constant in $G_0$:
indeed, replacing $G_0$ by $G_0+c$ changes the right-hand side of
\eqref{eq:res-equivariant-primitive} by
$$
                    \frac cM\sum_{g\in\mu_M}\chi_q(g)^{-1}=0.
$$
Equivalently, $G$ is the unique antiderivative of $W$ satisfying the second
identity in \eqref{eq:res-equivariant-primitive-properties}.    
\end{rem}

Since $\chi_q^M=1$ on $\mu_M$, the $M$-th power $G^M$ is deck-invariant.
\begin{defin}\label{eq:res-norm-definition}
The descent of the function $G^M$ to $C$ is called its \emph{norm}:
$$
 H\in\C(C),\qquad \pi^*H=G^M.
$$
\end{defin}

When $M=1$, put $G=H=G_0$.  For every $M$, the difference
$G_0-G$ is a constant, denoted by $c$; thus $c=0$ when $M=1$.

\begin{rem}[Geometric definition of the norm]
\label{rem:res-polygon-norm}
For $M>1$, the norm has the regular-polygon description from the
introduction.  Over a point $y\in C$ outside the branch points and poles,
let $g_1(y),\ldots,g_M(y)$ be the values of $G_0$ on the fiber.  They have
the form $c+\zeta^jG(\widetilde y)$, $0\le j<M$, where $\zeta$ is a
primitive $M$-th root of unity and $\widetilde y$ is any lift of $y$.
Thus they form a regular polygon with center
$c=M^{-1}\sum_jg_j(y)$; it collapses to its center when
$G(\widetilde y)=0$. Note that $c$, as a function of $y$, is constant, because its differential equals the average of $dG_0=W$, which is 0 due to the equivariance of $W$. Centering the vertices and multiplying gives
$$
 H(y)=(-1)^{M-1}\prod_{j=1}^M
       \left(g_j(y)-\frac1M\sum_{k=1}^M g_k(y)\right).
$$
Indeed, $\prod_{j=0}^{M-1}\zeta^j=(-1)^{M-1}$, so the right-hand side
is $G(\widetilde y)^M$.  \end{rem}

\subsection{The geometry of the norm} We describe zeros, poles, and local degrees of the norm $H$, in terms of the numerical data of the resolution.

\begin{lemma}
\label{lem:res-norm-local-degrees}
Assume that $C$ is clean.  Let $Q\subset C^\circ$ be the image of the
zeros of $G$ over $C^\circ$.  Partition the positive-defect punctures into
the set $\mathcal Z$ of those at whose lifts $G$ vanishes and the set
$\mathcal A$ of the remaining ones.  We claim every member of $\mathcal A$ is
unramified, and hence has $\alpha_p\in\Z_{>0}$
by~\eqref{eq:res-curve-inertia-formula}.
The zero and pole divisors of $H$ are
\begin{align}
 (H)_0&=M\sum_{x\in Q}x+\sum_{p\in\mathcal Z}(M\alpha_p)p,\notag\\
 (H)_\infty&=\sum_{\alpha_p<0}(-M\alpha_p)p,
                                                        \label{eq:res-norm-divisors}
\end{align}
the latter supported within the set of punctures $P$. The local degrees of $H$ are given by
\begin{equation}\label{eq:res-norm-local-degrees}
\begin{array}{c|c|c}
 \text{point of }C & \text{value of }H & \text{local degree of }H\\ \hline
 x\in Q &0&M\\
 p\in\mathcal Z&0&M\alpha_p\\
 p\in\mathcal A&\mbox{\rm in }\C^*&\alpha_p\\
 p\in P,\ \alpha_p<0&\infty&-M\alpha_p.
\end{array}
\end{equation}
At every other point the local degree is one.  In particular, the
ramification multiplicity at each point in the table is its displayed
local degree minus one.
\end{lemma}

\begin{proof}
Since $dG=W$, the function $G$ has no pole over
$C^\circ$, and its zeros there are simple.  At a lift $\widetilde p$
of a puncture $p$ with ramification index $e=e_p$, formula
\eqref{eq:res-curve-order-formula} gives
$$
 \ord_{\widetilde p}G=e\alpha_p\quad(\alpha_p<0),
 \qquad
 \ord_{\widetilde p}\bigl(G-G(\widetilde p)\bigr)
                       =e\alpha_p\quad(\alpha_p>0).
$$
If $\alpha_p>0$ and $e>1$, a nonidentity deck transformation fixes
$\widetilde p$ and acts on $G$ by a nontrivial scalar.  Therefore
$G(\widetilde p)=0$.  Thus every $p\in\mathcal A$ is unramified,
and its positive defect is an integer.  Equivariance shows that
the vanishing or nonvanishing of $G$ is the same at all lifts of $p$.

To make conclusions about $H$, we introduce local coordinates for each row of the table.  We shall use the
following elementary coordinate choice at a zero or pole of $G$.
Initially choose coordinates $z=\widetilde z^e$ for the covering, and write
$G=\widetilde z^kA(\widetilde z)$, where $k\in\Z\setminus\{0\}$ and
$A(0)\ne0$.  A generator of the group of deck transformations fixing
$\widetilde p$ acts by $\widetilde z\mapsto\zeta\widetilde z$, with $\zeta$
a primitive $e$-th root of unity.  Comparing the leading terms in the
equivariance identity for $G$ shows that
$A(\zeta\widetilde z)=A(\widetilde z)$.  Hence $A(\widetilde z)$ is a
holomorphic function of $z$.  Replacing $\widetilde z$ by
$\widetilde z A(\widetilde z)^{1/k}$, and $z$ by
its $e$-th power, gives coordinates in which
$$
                  z=\widetilde z^e,\qquad G=\widetilde z^k.
$$
The required roots of the unit $A$ exist on a sufficiently small
analytic neighborhood.

If $x\in Q$, the covering is unramified and $dG$ does not vanish at
a lift of $x$.  Thus $\widetilde z=G$ is a local coordinate there and also gives
a coordinate $z=\widetilde z$ downstairs.  In these coordinates
$$
                         G_0=c+\widetilde z,\qquad H=z^M.
$$
This proves the first row and its zero multiplicity.

If $p\in\mathcal Z$, put $k=e\alpha_p>0$.  The coordinate choice
above gives
$$
              z=\widetilde z^e,\qquad G_0=c+\widetilde z^k,
              \qquad H=z^{Mk/e}=z^{M\alpha_p}.
$$
This proves the second row and its zero multiplicity.

If $p\in\mathcal A$, put $b=G(\widetilde p)\ne0$ and
$k=\alpha_p>0$.  The covering is unramified, and
$G-b$ has a zero of order $k$.  After absorbing its leading unit into
a coordinate $z=\widetilde z$ we obtain
$$
 G_0=c+b+\widetilde z^k,\qquad
 H=(b+z^k)^M
   =b^M+Mb^{M-1}z^k+O(z^{2k}).
$$
Since $Mb^{M-1}\ne0$, the local degree is $k$, which proves the
third row.

If $\alpha_p<0$, put $k=e\alpha_p<0$.  The same coordinate
choice at the pole gives
$$
              z=\widetilde z^e,\qquad G_0=c+\widetilde z^k,
              \qquad H=z^{Mk/e}=z^{M\alpha_p}.
$$
Thus $H$ has a pole of order $-M\alpha_p$, and its reciprocal has a
zero of that order.  This proves the last row, including its local degree.

Every remaining point lies in $C^\circ$ and has $G(\widetilde x)=b\ne0$ with $dG(\widetilde x)=W(\widetilde x)\ne 0$.
Here $\widetilde z=G-b$ is a coordinate, and
$H=(b+z)^M$ has nonzero derivative at $z=0$.  These cases exhaust $C$
and prove both divisor formulas.  \end{proof}

\subsection{When can the curve residue be exact?}

\begin{lemma}\label{lem:res-exactness-rigidity}
If the residue $W$ on one connected covering component $\widetilde C$
is exact, then $C$ is clean. If moreover $B_{\mathrm{top}}(C)\ne0$, then some
puncture $p\notin C^\circ$ satisfies
\begin{equation}\label{eq:res-large-defect}
\alpha_p\in\Z,
 \qquad \alpha_p\ge2,
\end{equation}
so that in particular $\widetilde C$ is unramified over $p$. If $M>1$,  the
normalized equivariant antiderivative
\eqref{eq:res-equivariant-primitive} has a finite nonzero value above it.
\end{lemma}

The role of this lemma is to prove that the residue $W$ is not exact, thanks to the contradiction with the subsequent 
Proposition~\ref{prop:mmp-corner} asserting $\alpha_p\le 1$ in this setting. 

\begin{proof}
The differential of a meromorphic function cannot have a simple pole.
Formula~\eqref{eq:res-curve-order-formula} therefore excludes $\alpha_p=0$
at every puncture, so $C$ is clean.

Suppose first that $M>1$ and use the norm notation introduced above.
Equality of the degrees of the two divisors in
\eqref{eq:res-norm-divisors}, divided by $M$, gives
\begin{align}
 0\le|Q|&=-\sum_{\alpha_p<0}\alpha_p-
       \sum_{p\in\mathcal Z}\alpha_p                                      \notag\\
  &=\sum_{p\in\mathcal A}\alpha_p-
       \bigl(2g(C)-2+|P|\bigr)=\sum_{p\in\mathcal A}\alpha_p+e(C^\circ).                         \label{eq:res-zero-count}
\end{align}
The second equality is the adjunction identity
\eqref{eq:res-adjunction-identities}.

Assume that no member of $\mathcal A$ has defect at least $2$.  Every
member of $\mathcal A$ then has defect $1$, 
and formula
\eqref{eq:res-zero-count} becomes
$$0\le |Q|=e(C^\circ\cup\mathcal{A}),$$
meaning $C^\circ\cup\mathcal{A}$ is a projective line missing at most two points or a complete elliptic curve. Since $G$ must have a pole, and it can exist only above a puncture of $C^\circ\cup\mathcal{A}$, we further conclude $C^\circ\cup\mathcal{A}$ is a projective line with one or two punctures. For one puncture, Lemma~\ref{lem:res-equal-ratio-adjunction} gives
the sole defect $-1$; for two it gives a pair of opposite defects.  In
both cases the second identity in \eqref{eq:res-adjunction-identities}
gives $B_{\mathrm{top}}(C)=0$, a contradiction.

Now suppose that $M=1$.  Formula \eqref{eq:res-curve-order-formula} makes
every defect integral.  If no positive defect is at least $2$, the norm
$H=G_0:C\to\Pone$ has no ramification over a finite value: $W$ has no
zero on $C^\circ$, positive defects are $1$, and negative defects are
poles. If $b$ is the number of poles, all ramification lies over infinity,
and Riemann--Hurwitz gives, in terms of the Euler characteristic $e(C)$,
$$
 e(C)=2\deg H-(\deg H-b)=\deg H+b.
$$
Since $e(C)\le2$, $\deg H\ge1$, and $b\ge1$, one has
$e(C)=2$ and $\deg H=b=1$.  The only nonneutral defect is again
$-1$, so $B_{\mathrm{top}}(C)=0$, a contradiction.
\end{proof}

\subsection{A prescribed curve gives a Malgrange eigencycle}
\label{subsec:res-individual-localization}

\begin{theor}[Malgrange eigencycle on a prescribed equal-ratio curve]
\label{thm:res-individual-RV-curve}
Let $f\in\C[x_1,\ldots,x_n]$, $n\ge2$, and let $C$ be a component of a vertical curve $E_I$ with $\nu_i/N_i=q$ for $i\in I$, contained in $\rho^{-1}(0)$ on an algebraic embedded log
resolution of $(f=0)$.  Assume either that $C$ meets a further ratio-$q$
component of $(F=0)$, or that $C$ is clean and satisfies one of the following conditions:
\begin{enumerate}
\item $B_{\mathrm{top}}(C)\ne0$;
\item $g(C)=1$ and $\alpha_p=1$ at every puncture.
\end{enumerate}
Put
$$
                               \lambda=e^{-2\pi iq}.
$$
Then, for every analytic neighborhood $V$ of $C$ in $X^{\mathrm{an}}$ and
all sufficiently small $t\ne0$,
\begin{equation}\label{eq:res-localized-RV-class}
 \im\!\left[
 H_{n-1}\bigl(V\cap\rho^{-1}(M_{0,t});\C\bigr)
       \to H_{n-1}(M_{0,t};\C)
 \right]
 \cap\ker(\M_0-\lambda\cdot{\rm Id})\ne0.
\end{equation}
The class is supported arbitrarily close to the prescribed curve.
Moreover, on every connected component of its
$q$-effective cover there is a compact smooth loop $\gamma$ such that
$$
                         \int_\gamma W_{I,q}\ne0.
$$
\end{theor}

In the second case $B_{\mathrm{top}}(C)=0$. This case will be needed in
Section~\ref{sec:poles} for motivic contributions invisible under
topological specialization, including poles whose order drops.

\begin{proof}
If $C$ is not clean, choose $p\in P$ with $\alpha_p=0$.
By Lemma~\ref{lem:res-curve-inertia}, the $q$-effective cover is
unramified over $p$, and $W_{I,q}$ has a simple pole at every
lift $\widetilde p$.  A positively oriented small loop $\gamma$
around any such lift satisfies
$$
 \int_\gamma W_{I,q}=2\pi i\Res_{\widetilde p}W_{I,q}\ne0.
$$
Every connected covering component contains a lift of $p$, so
Proposition~\ref{prop:general-Fourier-Leray} proves the assertion.
We may therefore assume that $C$ is clean.

Suppose first that $B_{\mathrm{top}}(C)\ne0$.
Use the neighborhood
and resolution $\rho:X\to U$ fixed in Section~\ref{sec:covers}, and
apply Proposition~\ref{prop:mmp-corner} to this resolution and the divisor
$\Div(f)|_U$, whose pullback is $\Div(F)$. At each $p\in P$, exactly $n$ components of the reduced total transform
meet: the $n-1$ parents of $C$ and the unique component transverse to $C$.
Their $q$-defects are $(0,\ldots,0,\alpha_p)$.
Proposition~\ref{prop:mmp-corner} therefore gives $\alpha_p\le1$ for every
$p\in P$.
If the residue on any
connected component of the $q$-effective cover were exact,
Lemma~\ref{lem:res-exactness-rigidity} would produce a puncture
$p_0\in P$ with integral defect $\alpha_{p_0}\ge2$: contradiction.  Hence
the residue is nonexact on every component.  A meromorphic differential on
a punctured complete curve which is not exact has a nonzero period.
If a linear combination of loop periods is nonzero, one of those periods
is nonzero, so a single compact smooth loop suffices.
Proposition
\ref{prop:general-Fourier-Leray} now gives
\eqref{eq:res-localized-RV-class}.

Suppose next that $g(C)=1$ and every $\alpha_p=1$.  The latter implies that
the $q$-effective cover is unramified at every puncture and that the residue
extends to a nonzero holomorphic one-form on the compactification of each
connected covering component.  Such a form is not exact: if $W=dG$ for a
meromorphic function $G$, then a pole of $G$ would give a pole of $dG$;
hence $G$ would be holomorphic on a complete curve, and therefore constant.
Thus $W$ again has a nonzero period, and
Proposition~\ref{prop:general-Fourier-Leray} applies.
\end{proof}

\section{The corner inequality}
\label{sec:corner}

We prove the numerical inequality $(0,\ldots,0,a)\Rightarrow a\le1$
from Theorem~\ref{thmmp}. Two statements extracted from the proof of
\cite[Theorem~2.4]{NicaiseXu} provide a sequence of birational maps and
show that every coefficient-one boundary component meeting the contracted
fiber through the chosen point also contains that point. The remaining
argument finds a curve in this fiber and bounds two intersection numbers
on a surface. We recall the required definitions below and use the cited
proof for the construction of the birational maps.

Let $U$ be a smooth irreducible complex
variety of dimension $n$, let
$$
                       S=\sum_a m_aS_a
$$
be a nonzero effective integral Cartier divisor with distinct reduced prime
components $S_a$, and let
$$
                         h:X\to U
$$
be a projective embedded log resolution which is an isomorphism over
$U\setminus\Supp S$.  Write
\begin{align}
 h^*S&=\sum_a m_aS_{a,X}
          +\sum_{F\subset\Exc(h)}N_FF,
 \label{eq:mmp-resolution-N}\\
 K_{X/U}&=\sum_{F\subset\Exc(h)}(\nu_F-1)F.
 \label{eq:mmp-resolution-nu}
\end{align}
For a strict transform $S_{a,X}$ set $(N,\nu)=(m_a,1)$.  Fix
$q\in\Q_{>0}$ and define the $q$-defect of every component $G$ of
the reduced total transform by
\begin{equation}\label{eq:mmp-defect}
                         \alpha_q(G)=\nu_G-qN_G.
\end{equation}

\begin{utver}[Corner inequality]\label{prop:mmp-corner}
Assume $n\ge2$.  Let $p\in X$ lie on exactly $n$ irreducible components
$G_1,\ldots,G_n$ of
$(h^{-1}(\Supp S))_{\mathrm{red}}$.  These components are
smooth and meet transversely at $p$, because $h$ is an embedded log
resolution.  If, after reordering them,
\begin{equation}\label{eq:mmp-corner-data}
 \bigl(\alpha_q(G_1),\ldots,\alpha_q(G_n)\bigr)
                         =(0,\ldots,0,a),
\end{equation}
then $a\leq1$.
\end{utver}

No simple-normal-crossing hypothesis is imposed on $S$ or on
$S_{\mathrm{red}}$ in $U$.  The proof is given in
Subsection~\ref{subsec:mmp-corner-proof}; it uses the preparations in
Subsections~\ref{subsec:mmp-pairs}, \ref{subsec:mmp-capped}, and
\ref{subsec:mmp-local-index}.

The statement is local on $U$, so we may shrink $U$ to an affine
open neighborhood of $h(p)$ on which $S$ is the divisor of a regular
function, as in the setup of \cite[\S2]{NicaiseXu}.

\subsection{Birational preliminaries}
\label{subsec:mmp-pairs}

For a proper birational morphism $g:Y\to Z$, the exceptional locus
$\Exc(g)$ is the locus where it is not a local isomorphism. A
\emph{contraction} in what follows is such a projective morphism with
connected fibers. 

A \emph{boundary} $\Delta$ in a normal variety $Y$ is an effective $\Q$-divisor with coefficients
in $[0,1]$. A \emph{pair} $(Y,\Delta)$ consists of a normal variety and a
boundary for which $K_Y+\Delta$ is $\Q$-Cartier. A log resolution
$\mu:\overline Y\to Y$ of the pair is a proper birational morphism with
$\overline Y$ smooth, $\Exc(\mu)$ a divisor and
$\Exc(\mu)\cup\Supp(\mu^{-1}_*\Delta)$ a simple-normal-crossing divisor.
Here $\mu^{-1}_*\Delta$ is the strict transform. For each prime divisor
$F$ on such a resolution, including strict transforms, define its
\emph{log discrepancy} by
\begin{equation}\label{eq:mmp-discrepancies}
 K_{\overline Y}
   =\mu^*(K_Y+\Delta)
      +\sum_F\bigl(A(F;Y,\Delta)-1\bigr)F.
\end{equation}
A prime divisor over $Y$ means one on a normal variety with a proper
birational morphism to $Y$; its \emph{center} on $Y$ is the closure of its
image. To define $A$ for such a divisor, choose a log resolution dominating
that model and use \eqref{eq:mmp-discrepancies} for its strict transform.
This transform is nonempty: the morphism from the resolution to that normal
model is an isomorphism near the generic point of the divisor, even if the
divisor is exceptional over $Y$.
The result is independent of this choice. A boundary component of
coefficient $b$ has $A=1-b$.

The pair is \emph{log canonical} (lc) if it admits a log resolution on
which every exceptional divisor has $A\ge0$. This condition then holds
on every log resolution. A \emph{log canonical
center} is the center of a prime divisor over $Y$ with $A=0$.

The pair is \emph{divisorial log terminal} (dlt) if it
admits a log resolution on which every exceptional divisor has $A>0$.
The existence of one such resolution is the condition; further refinement of such a resolution (particularly, blowing up an intersection of coefficient-one components) can create a divisor with $A=0$. 

For a dlt pair, write $\lfloor\Delta\rfloor$ for the reduced union of its
coefficient-one components. Its lc centers are exactly the irreducible
components of intersections of these components. Such centers are normal,
and the pair is SNC at their generic points. Moreover, intersections of a
fixed collection of coefficient-one components have disjoint irreducible
connected components, see e.g. \cite[Theorem~4.16]{KollarSingularities}. Thus an lc
curve is smooth, although the ambient variety may be singular at some of
its points. A zero-dimensional lc center is an SNC point at which exactly
$\dim Y$ coefficient-one components meet.

For a dlt pair $(Y,\Delta)$ and a coefficient-one component $F$, the
\emph{adjunction theorem} supplies an effective boundary $\Delta_F$ such
that $(F,\Delta_F)$ is dlt and
$$
                (K_Y+\Delta)|_F\sim_{\Q}K_F+\Delta_F;
$$
see \cite[Definition--Lemma~3.4.1(1)]{BCHM}. The boundary supplied by
adjunction is called the \emph{different}; the displayed linear equivalence
alone does not define it. Restriction here means
restriction of a sufficiently divisible Cartier multiple of
$K_Y+\Delta$, and the symbol
$\sim_{\Q}$ means that suitable positive integer multiples are linearly
equivalent. Iterating
along an ordered collection of coefficient-one components cutting out a
center $W$ gives
\begin{equation}\label{eq:mmp-dlt-adjunction}
       (K_Y+\Delta)|_W\sim_{\Q}K_W+\Diff_W(\Delta).
\end{equation}
The effective boundary $\Diff_W(\Delta)$ is called the \emph{different}.
On a curve it is a sum of points with rational coefficients; at an SNC
point a further coefficient-one component transverse to $W$ contributes
one.

\subsection{The Nicaise--Xu toolkit}
\label{subsec:mmp-capped}

We use the following two steps of the proof of
\cite[Theorem~2.4]{NicaiseXu}.  They do not require the equal-weight
hypothesis of that theorem.  Put
$$
 \Delta_q:=\Exc(h)_{\mathrm{red}}
                +\sum_a\min\{qm_a,1\}S_{a,X},
$$
where $m_a$ and $S_{a,X}$ are as in \eqref{eq:mmp-resolution-N}.
Thus every exceptional component has coefficient one, and so does every
strict component of zero $q$-defect.

\begin{lemma}\label{lem:mmp-nx-sequence}
There is a finite sequence of birational maps over $U$
$$
 (X,\Delta_q)=(X_0,\Delta_0)\dashrightarrow\cdots
       \dashrightarrow(X_m,\Delta_m),\qquad h_j:X_j\to U,
$$
with the following properties.
\begin{enumerate}
\item Each $h_j$ is projective and birational, $X_j$ is normal and
$\Q$-factorial, $(X_j,\Delta_j)$ is dlt, and $\Delta_j$ is the
birational transform of $\Delta_q$. Every prime divisor on $X_{j+1}$ is
the birational transform of a prime divisor on $X_j$.
\item Each step has an associated projective birational morphism
$g_j:X_j\to Z_j$ over $U$, with $Z_j$ normal and with connected fibers.
The birational map $X_j\dashrightarrow X_{j+1}$ is an isomorphism on
$X_j\setminus\Exc(g_j)$.  The numerical classes over $U$ of
all complete curves contracted by $g_j$ are positive multiples of one
another\footnote{For clarity, the numerical class of a curve over $U$ records its degrees
against all Cartier divisors on the model; only curves mapped to a point
of $U$ are considered.  Thus the proportionality in the lemma means that,
for two curves $\Gamma$ and $\Gamma'$ contracted by the same $g_j$,
there is a number $b>0$ such that
$L\cdot\Gamma'=bL\cdot\Gamma$ for every $\Q$-Cartier divisor $L$ on $X_j$.
No choice of a homology or Chow realization is needed here.}, and every such curve $\Gamma$ satisfies
$$
                   (K_{X_j}+\Delta_j)\cdot\Gamma<0.
$$
\item On the final model,
$$
       K_{X_m}+\Delta_m-h_m^*(K_U+qS)\leq0.
$$
In particular, every exceptional component of $X$ with positive
$q$-defect disappears in this sequence.
\end{enumerate}
\end{lemma}
\begin{rem}
When we later apply this lemma to prove Proposition~\ref{prop:mmp-corner}, assuming towards contradiction that $a>1$, the last sentence of the lemma assures that some step of the sequence of maps changes a neighborhood of $p$; the proof uses the first such step and the properties of its source pair supplied by the lemma. \end{rem}

The proof of the lemma is part of the proof of \cite[Theorem~2.4]{NicaiseXu}.
The sequence and property~(1) are established before the Claim, on p.~5
of arXiv:1403.6792v2. Property~(2) appears in the proof of the Claim, around (2.5). Property~(3) is the assertion that $-D$ is effective in the
paragraph following that proof, on p.~6. These steps use no equal-weight
assumption.  The last sentence follows because the coefficient of
a surviving exceptional divisor $F$ in the left-hand side is
$\nu_F-qN_F$.

We shall use one more observation from the referred proof.

\begin{lemma}[The fiber through a point stratum]
\label{lem:mmp-nx-fiber}
Let $g_j$ be a contraction in Lemma~\ref{lem:mmp-nx-sequence}, and
let $p\in X_j$ be a zero-dimensional log canonical center of
$(X_j,\Delta_j)$.  Every log canonical center meeting
$g_j^{-1}(g_j(p))$ contains $p$.  In particular, every component of
$\lfloor\Delta_j\rfloor$ meeting this fiber passes through $p$, and
the fiber contains no other zero-dimensional log canonical center.
\end{lemma}

This is the argument following (2.5) in the proof of the Claim (the last
paragraph of p.~5 and the first paragraph of p.~6 of arXiv:1403.6792v2) in
\cite[Theorem~2.4]{NicaiseXu}.  More precisely, take the center denoted
$C_\ell$ therein to be the point $p$.  The uniqueness of the smallest
center meeting the fiber applies because the point $p$ is itself such
a center; the subsequent argument excludes all coefficient-one
components not containing it.  This step of the proof uses
no numerical equality of weights.

\begin{lemma}[Pairing a point with a curve]\label{lem:mmp-paired-curve}
Let $g_j:X_j\to Z_j$ be a contraction in
Lemma~\ref{lem:mmp-nx-sequence}, and let $p\in\Exc(g_j)$ be a
zero-dimensional lc center. Suppose that a coefficient-one boundary
component $D_0$ through $p$ satisfies $D_0\cdot\Gamma>0$ for the curves contracted
by $g_j$. There is a one-dimensional lc center $W\subset\Exc(g_j)$
not contained in $D_0$ such that $p$ is an irreducible component of
$D_0\cap W$. Moreover, $D_0\cap W$ is nonempty and irreducible and maps
onto $g_j(W)$.
\end{lemma}

These are Cases~1 and~2 in the proof of
\cite[Theorem~19]{dFKX}, applied to the contracted center $\{p\}$.
Case~1 supplies $W$ and its intersection with $D_0$; Case~2 supplies
irreducibility and the assertion about its image. The contractions here
are those of the actual MMP sequence constructed in
\cite[proof of Theorem~2.4]{NicaiseXu}, as required by that result.

\subsection{The numerical bound on a boundary curve}
\label{subsec:mmp-local-index}

Recall
that the local Cartier index of a $\Q$-Cartier divisor at a point is the
least positive integer whose multiple is Cartier near that point.
We use the following surface adjunction fact
\cite[Definition--Lemma~3.4.1(3)]{BCHM}: if $(Z,W)$ is a dlt surface pair
whose boundary is a single reduced irreducible curve, then $W$ is smooth
and the adjunction boundary $\Theta_0$ in
$(K_Z+W)|_W\sim_{\Q}K_W+\Theta_0$ has coefficient $1-1/r_x$ at $x\in W$,
where $r_x$ is the local Cartier index of $W$.  The hypothesis used in the cited
statement is equivalent, for a dlt pair, to having disjoint coefficient-one
boundary components. 

\begin{lemma}[One transverse boundary component]\label{lem:mmp-one-floor}
Let $(Y,\Delta)$ be a $\Q$-factorial dlt $\Q$-pair of dimension $n\ge2$
over $\C$, and let $W$ be a complete one-dimensional log canonical center.
Assume that $(K_Y+\Delta)\cdot W<0$, and that exactly one component of
$\lfloor\Delta\rfloor$ not containing $W$ meets $W$.  Assume that its
intersection with $W$ is the single point $p$, where $(Y,\Delta)$ is SNC.
If a component $E\supset W$ of $\lfloor\Delta\rfloor$ satisfies
$E\cdot W<0$, then
\begin{equation}\label{eq:mmp-length-one}
       \frac{-(K_Y+\Delta)\cdot W}{-E\cdot W}\leq1.
\end{equation}
\end{lemma}

\begin{proof}
At the generic point of $W$, the pair is SNC, and exactly $n-1$ floor
components contain $W$.  Besides $E$, write them as
$F_1,\ldots,F_{n-2}$.  Let $Z$ be the connected component of
$F_1\cap\cdots\cap F_{n-2}$ containing $W$; for $n=2$, put $Z:=Y$.
Intersections of a fixed collection of coefficient-one components have disjoint
irreducible normal connected components
\cite[Theorem~4.16]{KollarSingularities}.  Thus $Z$ is
a normal surface and, on a neighborhood of $W$, the restriction of the
$\Q$-Cartier divisor $E$ to $Z$ is precisely $W$.  Indeed, this restriction
has multiplicity one at the generic point of $W$, where the pair is SNC,
and all its other components are disjoint from $W$.  In particular,
$$
                         E\cdot W=W_Z^2.
$$

Iterated dlt adjunction \eqref{eq:mmp-dlt-adjunction} gives an effective
boundary $\Delta_Z$ such that
$(Z,\Delta_Z)$ is dlt and
$$
 (K_Y+\Delta)|_Z\sim_{\Q}K_Z+\Delta_Z,
 \qquad (K_Y+\Delta)\cdot W=(K_Z+\Delta_Z)\cdot W.
$$
The curve $W$ has coefficient one
in $\Delta_Z$ and is smooth.  Write its adjunction formula as
$$
 (K_Z+\Delta_Z)|_W\sim_{\Q}K_W+\Theta,
 \qquad \Theta=\sum_x\theta_x[x]\ge0,
$$
where $x$ runs through the points of $W$ and only finitely many
$\theta_x$ are nonzero. At $p$, the $n$ coefficient-one branches exhaust
the SNC boundary, so $\theta_p=1$. The inequality
$$
             \deg(K_W+\Theta)=(K_Y+\Delta)\cdot W<0
$$
then implies $W\simeq\Pone$ and
\begin{equation}\label{eq:mmp-adjunction-degree}
 \sum_{x\ne p}\theta_x<1,
 \qquad
 -(K_Y+\Delta)\cdot W=1-\sum_{x\ne p}\theta_x.
\end{equation}

The dlt surface $Z$ is $\Q$-factorial
\cite[Proposition~4.11]{KM}, so $(Z,W)$ is dlt after deleting the other
boundary components.  Apply the surface adjunction fact to $(Z,W)$.
Since $\Delta_Z-W$ is effective and has no component equal to $W$, its
restriction to $W$ is effective.  Thus, for $x\in W\setminus\{p\}$,
\begin{equation}\label{eq:mmp-different-index}
                         \theta_x\ge1-\frac1{r_x}.
\end{equation}
Here we use the additivity of adjunction under adding a $\Q$-Cartier
boundary: $\Theta=\Theta_0+(\Delta_Z-W)|_W$.

A point with $r_x>1$ contributes at least $1/2$ to the sum in
\eqref{eq:mmp-adjunction-degree}; hence there is at most one such point.
Let its index be $r$, or put $r:=1$ if there is none.  The divisor $rW$
is Cartier on a neighborhood of $W$ in $Z$, since it is Cartier at that
point and the local index is one everywhere else on $W$ (also at the SNC
point $p$).  Consequently,
\begin{equation}\label{eq:mmp-E-index}
              -E\cdot W=-W_Z^2\in\frac1r\Z_{>0},
              \qquad -E\cdot W\ge\frac1r.
\end{equation}
On the other hand, \eqref{eq:mmp-adjunction-degree} and
\eqref{eq:mmp-different-index} give
$$
                   0<-(K_Y+\Delta)\cdot W\le\frac1r.
$$
Dividing by \eqref{eq:mmp-E-index} proves the lemma.
\end{proof}

\begin{rem}
For $n=2$, the surface $Z$ is $Y$ itself and $E=W$; the lemma says
$-(K_Y+\Delta)\cdot W\le-W^2$.  In higher dimensions, the reduction to $Z$ preserves both intersection numbers. \end{rem}

\subsection{Proof of the corner inequality}
\label{subsec:mmp-corner-proof}

\begin{proof}[Proof of Proposition~\ref{prop:mmp-corner}]
Assume, for a contradiction, that $a>1$, and put $E:=G_n$ for the
component with defect $a$. A strict
component has defect $1-qm_b<1$, so $E$ is exceptional. In the boundary
$\Delta_q$ of Lemma~\ref{lem:mmp-nx-sequence}, both $E$ and the $n-1$
zero-defect components have coefficient one. Their SNC intersection $p$
is therefore a zero-dimensional lc center.

\smallskip\noindent\emph{The first contraction meeting $p$.}
By Lemma~\ref{lem:mmp-nx-sequence}(3), the positive-defect divisor $E$
must disappear. Thus there is a first step of that sequence which is not
an isomorphism on a neighborhood of $p$. The preceding steps identify
such a neighborhood with its image in $X_j$, and the same labeled SNC
corner remains there. We keep the notation $G_1,\ldots,G_n$ for their
transforms on $X_j$. Let $g:=g_j:X_j\to Z_j$ be the associated
contraction. In what follows, exceptional components and strict transforms
are taken with respect to $h_j:X_j\to U$. Then $p\in\Exc(g)$, and the fiber
$g^{-1}(g(p))$ contains a complete irreducible curve $\Gamma$
through $p$. Indeed, a proper birational morphism with connected fibers
onto a normal variety is an isomorphism wherever it is quasi-finite;
at $p$ the fiber is therefore positive-dimensional, and its projectivity
supplies such a curve.

By \eqref{eq:mmp-resolution-N}--\eqref{eq:mmp-resolution-nu}, an
exceptional divisor $F$ has coefficients $N_F$ in $h^*S$ and $\nu_F-1$
in $K_X-h^*K_U$, whereas a strict transform $S_{b,X}$ has coefficients
$m_b$ and $0$. These coefficients are unchanged for a surviving divisor,
since the birational maps are isomorphisms near its generic point.
Subtracting $q h_j^*S$ from $K_{X_j}-h_j^*K_U$ and adding $\Delta_j$,
whose coefficients are $1$ and $\min\{qm_b,1\}$ respectively by
Lemma~\ref{lem:mmp-nx-sequence}(1), gives
\begin{align}
 &K_{X_j}+\Delta_j-h_j^*(K_U+qS)\notag\\
 &\qquad=\sum_{F\subset\Exc(h_j)}(\nu_F-qN_F)F
       +\sum_b\bigl(\min\{qm_b,1\}-qm_b\bigr)S_{b,j}.
       \label{eq:mmp-Aq}
\end{align}
Only prime divisors present on $X_j$ occur in these sums, and $S_{b,j}$
denotes the strict transform of $S_b\subset U$. The subtracted pullback
$h_j^*(K_U+qS)$ has degree zero on every curve mapped to a point of $U$.

By Lemma~\ref{lem:mmp-nx-fiber}, every coefficient-one component meeting
this fiber contains $p$. An exceptional component not among
$G_1,\ldots,G_n$ therefore misses this fiber. A strict-transform term in
\eqref{eq:mmp-Aq} either has coefficient zero, if $qm_b\le1$, or is a
coefficient-one component missing this fiber, if $qm_b>1$. Indeed, in the
latter case it cannot be among $G_1,\ldots,G_n$: the component $G_n=E$
is exceptional, and a strict transform among $G_1,\ldots,G_{n-1}$ has
zero defect, hence $qm_b=1$. The components $G_1,\ldots,G_{n-1}$ have
zero defect. Consequently,
\begin{equation}\label{eq:mmp-negative-E}
  0>(K_{X_j}+\Delta_j)\cdot\Gamma
        =aE\cdot\Gamma.
\end{equation}
It follows that $E\cdot\Gamma<0$, and hence $\Gamma\subset E$.

\smallskip\noindent\emph{The contracted boundary curve.}
At least one of $G_1,\ldots,G_{n-1}$, denoted by $D_0$, satisfies
$D_0\cdot\Gamma>0$. Otherwise $\Gamma$ would be contained in each of
them: an effective $\Q$-Cartier divisor through $p$ has strictly
positive degree on $\Gamma$ unless it contains $\Gamma$, as one sees
by restricting an effective Cartier multiple to the normalization of
$\Gamma$. Together with $\Gamma\subset E$, this would put a curve
through $p$ inside the transverse intersection of $n$ divisors, which is
zero-dimensional near $p$.

The proportionality of the contracted curve classes in
Lemma~\ref{lem:mmp-nx-sequence} shows that $D_0$ has positive degree on
every curve contracted by $g$. Lemma~\ref{lem:mmp-paired-curve} supplies
a one-dimensional lc center $W\subset\Exc(g)$, not contained in $D_0$,
for which $p$ is an irreducible component of $D_0\cap W$. The same lemma
says that this intersection is irreducible and maps onto $g(W)$.
Therefore
\begin{equation}\label{eq:mmp-paired-W}
                     D_0\cap W=\{p\},\qquad g(W)=\{g(p)\}.
\end{equation}
In particular $W$ is complete, since it is a closed curve in the
projective fiber through $p$. 

At the SNC point $p$, the curve $W$ is the intersection of the $n-1$
components among $G_1,\ldots,G_n$ other than $D_0$. Thus $W\subset E$, and these are
exactly the coefficient-one components containing $W$. A further such
component meeting $W$ away from $p$ would give another zero-dimensional
lc center in this fiber, contrary to Lemma~\ref{lem:mmp-nx-fiber}. At $p$
there is no additional component by the SNC condition. Hence $D_0$ is
the unique transverse coefficient-one component meeting $W$, at the
single point $p$.

Applying \eqref{eq:mmp-Aq} on $W$ gives, by the same exclusion of the
other terms,
$$
  0<-(K_{X_j}+\Delta_j)\cdot W=a(-E\cdot W).
$$
Thus $E\cdot W<0$, and all the hypotheses of
Lemma~\ref{lem:mmp-one-floor} hold. It follows that
$$
        a=\frac{-(K_{X_j}+\Delta_j)\cdot W}{-E\cdot W}\le1,
$$
contrary to $a>1$.
\end{proof}

\begin{rem}
When $n=2$, the curve $\Gamma\subset E$ is $E$ itself, and the paired
curve $W$ is also $E$. The numerical estimate is then the surface
adjunction inequality $-(K_{X_j}+\Delta_j)\cdot E\le-E^2$.
\end{rem}

\section{High-order poles and global monodromy}
\label{sec:poles}

We extend our results from topological to motivic poles, using Rodrigues's geometric determination of
poles of the two highest possible orders \cite[Theorem~2.8.1 and \S3.2]{Rodrigues}.

\subsection{Motivic zeta functions and their poles}
\label{subsec:motivic-formulas}

The definition of the topological zeta-function in terms of a resolution \eqref{eq:topological-resolution-formula}
can be refined by replacing the Euler characteristics of the strata with stronger additive invariants.
Put
$$
       \mathcal M_{\C}=K_0(\Var_{\C})[\Lef^{-1}],
       \qquad \Lef=[\A^1].
$$
For $m\ge1$, let $\mathcal X_m$ be the space of $m$-jets $\varphi$ on
$U$ with $\varphi(0)=0$ and $\ord_t f(\varphi)=m$.  
\begin{defin}

The naive motivic zeta function is
$$
 Z^{\mathrm{naive}}_f(T)
   =\sum_{m\ge1}[\mathcal X_m]\Lef^{-mn}T^m.
$$    
\end{defin}
The change-of-variables formula gives
\begin{equation}\label{eq:naive-resolution-formula}
 Z^{\mathrm{naive}}_f(T)
  =\sum_{\varnothing\ne I\subset J}
   [E^\circ_I\cap\rho^{-1}(0)](\Lef-1)^{|I|}
   \prod_{i\in I}
    \frac{\Lef^{-\nu_i}T^{N_i}}
         {1-\Lef^{-\nu_i}T^{N_i}}.
\end{equation}
It lies in the subring of $\mathcal M_{\C}[[T]]$ generated by the displayed
rational factors; the right side is independent of the chosen resolution.
These formulas, including their local versions, are due to Denef and
Loeser \cite{DenefLoeser1992,DenefLoeser1998,DenefLoeser1999}.

Applying the compactly supported Hodge--Deligne homomorphism
$$
       E_c:K_0(\Var_{\C})\to\Z[u,v]
$$
to
\eqref{eq:naive-resolution-formula} gives the Hodge zeta
function
\begin{equation}\label{eq:Hodge-resolution-formula}
 Z^{\mathrm {Hdg}}_f(T)
  =\sum_{\varnothing\ne I\subset J}
   E_c(E^\circ_I\cap\rho^{-1}(0);u,v)(uv-1)^{|I|}
   \prod_{i\in I}
    \frac{(uv)^{-\nu_i}T^{N_i}}
         {1-(uv)^{-\nu_i}T^{N_i}}.
\end{equation} 

The Euler-characteristic limit recovers
\eqref{eq:topological-resolution-formula}: after $T=\Lef^{-s}$,
$$
 \lim_{\Lef\to1}
  (\Lef-1)
  \frac{\Lef^{-\nu_i-sN_i}}
       {1-\Lef^{-\nu_i-sN_i}}
     =\frac1{\nu_i+sN_i}.
$$

These three zeta functions have different coefficient rings and therefore support
different notions of pole.  For the rational functions
\eqref{eq:Hodge-resolution-formula} and
\eqref{eq:topological-resolution-formula}, one may use ordinary pole order.
For \eqref{eq:naive-resolution-formula}, a displayed
denominator in the motivic resolution formula does not by itself define a
pole, because the Grothendieck ring of varieties is not a domain.  To make up for it, adjoin compatible
rational powers of $\Lef$:
$$
 \mathcal M_{\C}^{\Q}
   =\varinjlim_{d\ge1}
      \mathcal M_{\C}[\Lef^{1/d}]
       \big/\bigl((\Lef^{1/d})^d-\Lef\bigr),
$$
where the positive integers are ordered by divisibility and, for $d\mid e$,
the transition map sends $\Lef^{1/d}$ to $(\Lef^{1/e})^{e/d}$. Put
\begin{equation}\label{eq:poles-RV-ring}
 \mathcal A_{\mathrm{RV}}
  =\mathcal M_{\C}^{\Q}
   \left[(1-\Lef^a)^{-1},m^{-1}
     \ \middle|\ a\in\Q\setminus\{0\},\;
                    m\in\Z\setminus\{0\}\right].
\end{equation}
We will work in the localization of $\mathcal A_{\mathrm{RV}}[T]$
obtained by inverting all factors
$$
        1-\Lef^{-a}T^b,
        \qquad a\in\Q,\quad b\in\Z_{>0}.
$$
Each polynomial $1-\Lef^{-a}T^b$ has invertible leading coefficient and hence is not a zero
divisor, even though the coefficient ring may have zero divisors.
For $0\ne H(T)\in\mathcal A_{\mathrm{RV}}[T]$ and $q\in\Q$, let
\begin{equation}\label{eq:poles-RV-valuation}
 v_q(H)=\max\{r\ge0:(T-\Lef^q)^r\text{ divides }H(T)
                   \text{ in }\mathcal A_{\mathrm{RV}}[T]\}.
\end{equation}
For a nonzero element $Z(T)$ of this localization, choose a representation
$$
               Z(T)=\frac{P(T)}{Q(T)},
$$
where $P\in\mathcal A_{\mathrm{RV}}[T]$ is nonzero and $Q$ is a finite
product of polynomials $1-\Lef^{-a}T^b$ with $a\in\Q$ and $b\in\Z_{>0}$.
Its Rodrigues--Veys order at $s=-q$ is
\begin{equation}\label{eq:poles-RV-order}
            \ord^{\mathrm{RV}}_{s=-q}Z=v_q(Q)-v_q(P).
\end{equation}
The divisibility result of Rodrigues--Veys says that this integer is
independent of the chosen representation; see
\cite[\S4]{RodriguesVeys} and \cite[\S3.1]{Rodrigues}.  \begin{defin}
We call $-q$
an \emph{RV pole of exact order $r$} of $Z(T)$ when
\eqref{eq:poles-RV-order} equals $r>0$. 
\end{defin}  The zero element is declared to
have no pole.
Note that the RV order does not test the vanishing at other points $T=\xi\Lef^q$, $\xi\ne1$ a root of unity.  
The Hodge pole
used below is defined in the same principal way: in the field
$$
        \bigcup_{d\ge1}\Q(u^{1/d},v^{1/d}),
$$
substitute $T=(uv)^q(1-\tau)$ in $Z^{\mathrm Hdg}_f(T)$ and write its
Laurent expansion as
$$
 \sum_{k\ge k_0}a_k\tau^k,\qquad a_{k_0}\ne0.
$$
Its \emph{principal Hodge order} at $-q$ is $-k_0$, the negative
of the smallest exponent with a nonzero coefficient. It is a pole when
$k_0<0$; the zero function has no pole.

The contribution of a 0-dimensional stratum of the exceptional divisor is given by the following straightforward fact; an extension to 1-dimensional strata occupies the next subsection.

\begin{exa}\label{cor:poles-maximal-motivic}
Let $x\in f^{-1}(0)$ and $q>0$. A principal Rodrigues--Veys pole of order
$n$ at $-q$, a principal Hodge pole of order $n$ there, and a
topological pole of order $n$ there are all equivalent to the existence of a point $p\in\rho^{-1}(x)$ at which $n$ components $E_i$ with $\nu_i/N_i=q$ meet. (These are components of the total transform and can include strict transforms.)
\end{exa}

\subsection{Individual curve contributions}
\label{subsec:res-motivic-bracket}

We group the resolution terms on a curve as in
\eqref{eq:topological-stratum-package}. In the naive motivic formula
\eqref{eq:naive-resolution-formula}, use the principal Rodrigues--Veys specialization
$T=\Lef^q(1-\tau)$.  For a clean vertical  component $C$ of an equal-ratio curve $E_I$,
grouping its open stratum with all its puncture strata gives the coefficient
\begin{equation}\label{eq:res-individual-RV-coefficient}
 \frac{(\Lef-1)^{n-1}}{\prod_{i=1}^{n-1}N_i}
                    B_{\mathrm{mot}}(C)
\end{equation}
of $\tau^{-(n-1)}$, where
\begin{equation}\label{eq:res-Bmot}
 B_{\mathrm{mot}}(C)
   :=[C^\circ]+\sum_{p\in P}
          \frac{\Lef-1}{\Lef^{\alpha_p}-1}
\end{equation}
in the Rodrigues--Veys coefficient ring.  Indeed, the factor of each $E_i,\, i\in I$, has
leading term $(N_i\tau)^{-1}$, while a transverse factor is regular at
$T=\Lef^q$ and has value
\begin{equation}\label{eq:res-transverse-RV-value}
 \frac{(\Lef-1)\Lef^{-\alpha_p}}
      {1-\Lef^{-\alpha_p}}
       =\frac{\Lef-1}{\Lef^{\alpha_p}-1}.
\end{equation}
The prefactor in \eqref{eq:res-individual-RV-coefficient} is a unit in
the ring $\mathcal A_{\mathrm{RV}}$ in \eqref{eq:poles-RV-ring}.
\begin{defin}
We call the vertical curve $C$ \emph{RV-contributing} if
$B_{\mathrm{mot}}(C)\ne0$.      
\end{defin}
Under the bound $\alpha_p\le1$, supplied by Proposition~\ref{prop:mmp-corner}, the difference between
this contribution and the topological contribution \eqref{eq:res-Btop} is explicit as follows.

\begin{utver}[Topological--motivic dichotomy]
\label{prop:res-motivic-dichotomy}
Let $C$ be a clean vertical component of a curve $E_I$ with equal ratio $q=\nu_i/N_i,\,i\in I$, and let $P:=C\setminus E^\circ_I$ be the set of its punctures. If the defect
$\alpha_p\le1$ for every $p\in P$, then:
\begin{enumerate}
\item If $B_{\mathrm{top}}(C)\ne0$, then $C\simeq\Pone$. Moreover, after
filling the neutral punctures (those with $\alpha_p=1$), at least three punctures remain and every
remaining defect lies in
$$
                              (-1,1)\setminus\{0\}.
$$
In particular, every nonneutral defect is nonintegral, $q\notin\Z$,
and every connected component of the $q$-effective cover has degree $M>1$.
\item If $B_{\mathrm{top}}(C)=0$, then exactly one of the following holds:
\begin{enumerate}
\item The curve $C$ has genus one and every defect is $1$; in this case
$$
                    B_{\mathrm{mot}}(C)=[C]\ne0;
$$
\item The curve $C$ is rational and, after filling the neutral punctures, either one puncture of defect $-1$ or two punctures of opposite defects within $(-1,1)$ remain.
In this case $B_{\mathrm{mot}}(C)=0$.
\end{enumerate}
\end{enumerate}
Consequently,
\begin{equation}\label{eq:res-good-curve-dichotomy}
 B_{\mathrm{mot}}(C)\ne0
 \quad\Longrightarrow\quad
 \left\{
 \begin{array}{l}
 B_{\mathrm{top}}(C)\ne0,\quad\text{or}\\
 g(C)=1\text{ and }\alpha_p=1\text{ for every }p.
 \end{array}\right.
\end{equation}
\end{utver}

\begin{proof}
Adjunction and the bound $\alpha_p\le1$ give
\begin{equation}\label{eq:res-deficit-sum}
                       \sum_{p\in P}(1-\alpha_p)=2-2g(C),
\end{equation}
so $g(C)\le1$.  If $g(C)=1$, every summand is zero.  Hence all defects are
$1$, $B_{\mathrm{top}}(C)=0$, and filling all neutral points gives
$B_{\mathrm{mot}}(C)=[C]$.  Its Hodge specialization is nonzero, so this
class remains nonzero in the Rodrigues--Veys coefficient ring.

Suppose $C\simeq\Pone$.  First assume
$B_{\mathrm{top}}(C)\ne0$ and fill every neutral puncture.  If $m$ points
remain, \eqref{eq:res-deficit-sum} gives
$$
                          \sum_{j=1}^m(1-\alpha_j)=2,
                          \qquad \alpha_j<1.
$$
The cases $m=1$ and $m=2$ give respectively the defect $-1$ and a pair
of opposite defects, and both have zero topological contribution.  Thus $m\ge3$.
For every $j$,
$$
 \alpha_j=m-2-\sum_{i\ne j}\alpha_i>-1,
$$
because each of the other $m-1$ defects is strictly smaller than $1$.
Cleanliness excludes zero.  If $q$ were integral, all
$\alpha_p=\nu_p-qN_p$ would be integral, which is impossible.  If $M=1$,
Lemma~\ref{lem:res-curve-inertia} would likewise make every defect
integral.  This proves (1).

It remains to assume $B_{\mathrm{top}}(C)=0$.  Not all defects can be
positive: then every contribution $1/\alpha_p-\alpha_p$ is nonnegative,
so $B_{\mathrm{top}}(C)=0$ would force every $\alpha_p=1$, contradicting
the deficit sum $\sum_p(1-\alpha_p)=2$ in
\eqref{eq:res-deficit-sum}.
There is exactly one negative defect, say $-b$.  Indeed, each negative
defect contributes more than $1$ to the left side of
\eqref{eq:res-deficit-sum}.  The same formula gives $0<b\le1$.
If $b=1$, the deficit sum leaves only neutral positive punctures, so the
nonneutral list is $\{-1\}$.  Assume henceforth that $0<b<1$.
Write every positive nonneutral defect as $1-x_i$, with $0<x_i<1$; then
$$
                              \sum_i x_i=1-b.
$$
For
$$
                              \psi(x)=\frac1{1-x}-(1-x)
$$
one has
\begin{equation}\label{eq:res-strict-superadditivity}
 \psi(x+y)-\psi(x)-\psi(y)
 =\frac{xy(2-x-y)}{(1-x)(1-y)(1-x-y)}>0
\end{equation}
whenever $x,y>0$ and $x+y<1$. Repeated application of
\eqref{eq:res-strict-superadditivity} gives
$$
 \sum_i\psi(x_i)\le\psi\Bigl(\sum_i x_i\Bigr)=\frac1b-b,
$$
with strict inequality if at least two $x_i$ occur. However,
$$
 B_{\mathrm{top}}(C)=b-\frac1b+\sum_i\psi(x_i)=0,
$$
so equality holds. Since $\sum_i x_i=1-b>0$, exactly one $x_i$ occurs;
it equals $1-b$, giving the positive defect $b$. Together with the case
$b=1$, this gives the lists $\{-1\}$ and $\{b,-b\}$, besides neutral defects.

The identity
$$
 \frac{\Lef-1}{\Lef^b-1}
       +\frac{\Lef-1}{\Lef^{-b}-1}=-(\Lef-1)
$$
together with $[\Pone\setminus P]=\Lef+1-|P|$ gives
$B_{\mathrm{mot}}(C)=0$ for the list $\{b,-b\}$ with any number of neutral
points.  For the one-point list $\{-1\}$ one similarly has
$\Lef+(\Lef-1)/(\Lef^{-1}-1)=0$.  This proves (2).
\end{proof}

\begin{sledst}\label{cor:poles-individual-motivic}
A clean vertical equal-ratio curve $C$ with $B_{\mathrm{mot}}(C)\ne0$
satisfies the hypotheses and thus the conclusion of
Theorem~\ref{thm:res-individual-RV-curve}. In particular, the same holds if its Hodge specialization is nonzero.
\end{sledst}

\begin{proof}
Proposition~\ref{prop:mmp-corner} gives $\alpha_p\le1$ at each puncture.
Proposition~\ref{prop:res-motivic-dichotomy} then gives one of the two
clean-curve alternatives of Theorem~\ref{thm:res-individual-RV-curve}.
\end{proof}

\begin{exa}
The elliptic alternative is genuinely motivic and can occur at $q=1$.
It has no plane analogue, because every curve contracted to a point of the
smooth plane is rational.
Let
\begin{equation}\label{eq:res-fermat-cubic-example}
                              f=x^3+y^3+z^3.
\end{equation}
Blowing up the origin gives an exceptional plane with data $(N,\nu)=(3,3)$.
The strict transform has data $(1,1)$, and their intersection is the smooth
plane cubic $C$.  Thus $q=1$, $C$ is elliptic, and $P=\varnothing$.  Hence
$$
                         B_{\mathrm{top}}(C)=\chi(C)=0,
                         \qquad B_{\mathrm{mot}}(C)=[C]\ne0.
$$
Its order-two principal coefficient is
$$
                       \frac{(\Lef-1)^2}{3}[C],
$$
whose Hodge specialization is nonzero.
The point $s=-1$ is therefore an RV pole of exact order two of the
motivic and Hodge zeta functions, although the corresponding double
topological coefficient vanishes. In fact
$Z_{\mathrm{top}}(f;s)=(1+s)^{-1}$, so the topological pole is simple.  A Malgrange eigencycle near $C$ gives a nonzero localized
eigenvalue-$1$ class in $H_2(M_0;\C)$.  \end{exa}

\subsection{The global motivic theorem}

We formulate the part of Rodrigues's theorem concluding the proof of our main result. A
\emph{$q$-component} is a component $E_i$ with
$\nu_i/N_i=q$.  A \emph{$q$-curve} is an irreducible component of
an intersection of $n-1$ such components which meets
$\rho^{-1}(0)$.  It is \emph{vertical} if it is contained in
$\rho^{-1}(0)$, and \emph{horizontal} otherwise.  \begin{defin}
A \emph{special
point} is a point of $\rho^{-1}(0)$ lying on exactly $n-1$
$q$-components such that their intersection curve is horizontal.
A further component through the point is allowed, but must have a different ratio.    
\end{defin}  

\begin{theor}[\cite{Rodrigues}, Theorem~2.8.1 and \S3.2]\label{thm:poles-Rodrigues}
Let $n\ge3$, and let $\rho$ be a projective algebraic embedded resolution.
A Hodge pole $-q$ of order at least $n-1$ exists if
and only if at least one of the following configurations occurs:
\begin{enumerate}
\item a special point;
\item a nonrational vertical $q$-curve;
\item a rational vertical $q$-curve having at least three non-neutral punctures; \item the incidence graph of rational vertical $q$-curves contains a
      cycle.
\end{enumerate}
The same geometric criterion holds for the naive motivic zeta function,
with pole understood in the Rodrigues--Veys sense
\eqref{eq:poles-RV-order}.
\end{theor}

We now prove the monodromy eigenvalue version of our main result, leaving the Bernstein-Sato version for the next section (see
Theorem~\ref{thm:bs-high-order}). We prefer to separate the two for the reader's convenience, and because the eigenvalue version specifies the dimension and support of the explicitly constructed eigencycle, which is not straightforward to see from the Bernstein-Sato perspective.

\begin{theor}[Poles of exact order $n-1$]
\label{thm:poles-global-RV}
Let $f\in\C[x_1,\ldots,x_n]$, $n\ge2$, satisfy $f(0)=0$; no
isolatedness or reducedness assumption is made.  If $-q$ is an RV pole
of exact order $n-1$ of the local naive motivic zeta function at the
origin, then, with $\lambda=e^{-2\pi iq}$, one of the following holds.
\begin{enumerate}
\item There is a clean vertical equal-ratio curve $C$ satisfying one
of the two clean-curve conditions in Theorem~\ref{thm:res-individual-RV-curve},
such that, for every analytic neighborhood $V$ of $C$ in $X$ and all
sufficiently small $t\ne0$ relative to $V$,
\begin{equation}\label{eq:poles-global-vertical}
 \im\!\left[
 H_{n-1}(V\cap\rho^{-1}(M_{0,t});\C)\to H_{n-1}(M_0;\C)
 \right]\cap\ker(\M_0-\lambda)\ne0.
\end{equation}
\item There is a horizontal equal-ratio curve, and the nearby-point
conclusion \eqref{eq:poles-horizontal-eigenspace}, with its local support
statement, holds along a curve through the origin.
\end{enumerate}
In particular, the same alternatives and support statements
hold if $-q$ is a Hodge or topological
pole of exact order $n-1$.
\end{theor}

\begin{proof}
We prove the motivic assertion here; the Hodge and topological corollaries are straightforward. 
An exact order-$(n-1)$ RV pole admits no point over the origin at which $n$
ratio-$q$ components meet: such a point would give order $n$ by
Corollary~\ref{cor:poles-maximal-motivic}.

Apply Theorem~\ref{thm:poles-Rodrigues} for $n\ge3$. For $n=2$, the same
criterion applies by \cite[Remark~2.8.2(i) and \S3.2]{Rodrigues}: its exceptional
germ $x_1^Nx_2^N$ has a pole of order two, already excluded.
A special point gives the
horizontal alternative by Proposition~\ref{prop:poles-horizontal}.
The cycle alternative is impossible, as the intersection point of two participating curves would contribute the pole of order $n\ne n-1$.

Suppose next that Rodrigues supplies a nonrational vertical $q$-curve
$C$.  The absence of an $n$-fold ratio-$q$ point makes $C$
clean.  Proposition~\ref{prop:mmp-corner} gives
$\alpha_p\le1$ at every puncture, while
\eqref{eq:res-adjunction-identities} gives
$$
                 2g(C)-2=\sum_p(\alpha_p-1)\le0.
$$
Since $C$ is nonrational, it follows that $g(C)=1$ and every
$\alpha_p=1$.  The individual-curve Theorem
\ref{thm:res-individual-RV-curve} proves \eqref{eq:poles-global-vertical}.

It remains to consider a rational vertical curve having at least three
nonneutral incidences.  Again it is clean, and
Proposition~\ref{prop:mmp-corner} gives $\alpha_p\le1$ at every puncture.
If $B_{\mathrm{top}}(C)=0$,
Proposition~\ref{prop:res-motivic-dichotomy} says that a rational curve has
at most two nonneutral defects, a contradiction.  Hence
$B_{\mathrm{top}}(C)\ne0$. Theorem
\ref{thm:res-individual-RV-curve} again gives
\eqref{eq:poles-global-vertical}.  \end{proof}

\section{From residue periods to Bernstein--Sato roots}
\label{sec:bernstein-sato}

We now strengthen the eigenvalue assertion for poles of exact order
$n-1$ to an assertion about the pole itself. The additional estimate concerns every holomorphic
numerator on the Leray tubes.  The respective assertion for the maximal-order root is already
proved in \cite[Theorem~3.5 and \S3.6]{NicaiseXu}.

\subsection{Relation to Bernstein--Sato roots}

For $x\in f^{-1}(0)$, the local Bernstein--Sato polynomial $b_{f,x}(s)$ is the monic generator of the nonzero ideal of polynomials $b(s)\in\C[s]$
for which there is an identity
\begin{equation}\label{eq:bs-identity}
 \exists P \;:\; b(s)f^s=P(s)f^{s+1}.
\end{equation}
Here $P(s)$ is a differential operator with holomorphic germ
coefficients at $x$, polynomial in $s$, and $f^s$ is a formal symbol
with $\partial_{x_i}f^s=s(\partial_{x_i}f)f^{s-1}$.  Existence is the
Bernstein--Sato theorem; the relation with periods is developed in
\cite{MalgrangeBernstein}, and treated for the local Milnor fibration in \cite[\S6]{Blanco}.  As before,
$\omega=dx_1\wedge\cdots\wedge dx_n$.

\begin{theor}[Poles of exact order $n-1$ give roots]
\label{thm:bs-high-order}
Let $f\in\C[x_1,\ldots,x_n]$ be nonconstant, with $f(0)=0$, and let
$n\ge2$ and $q\in\Q_{>0}$.  If $-q$ is an RV pole of exact order
$n-1$ of $Z^{\mathrm{naive}}_f$ (in particular, if $-q$ is a topological zeta-function pole of exact order $n-1$),
then
$$
                         b_{f,0}(-q)=0.
$$
\end{theor}

The root is
at the origin even when the proof in Section~\ref{sec:poles} supplies a Malgrange eigencycle near a
horizontal curve in nearby Milnor fibers.  We first record the calculation with the Leibniz rule for the compact-period argument of \cite[Proposition~3.3]{MalgrangeBernstein}.

\begin{lemma}\label{lem:bs-leibniz}
In a Bernstein identity~\eqref{eq:bs-identity}, let $L\ge0$ bound the
degrees in $s$ of the coefficients of $P(s)$, and write
$P(s)=\sum_{\ell=0}^L s^\ell P_\ell$, where the $P_\ell$ are differential
operators with holomorphic germ coefficients.  Let $r\ge1$ bound their
differential orders.  There are holomorphic germs $h_\ell$ and holomorphic
$(n-1)$-forms $\eta_j(s)$, polynomial in $s$, such that
\begin{equation}\label{eq:bs-green-identity}
 \begin{aligned}
 b(s)f^s\omega&=(P(s)f^{s+1})\omega
                 =f^{s+1}h_s\omega+d\Xi_s,\\
 h_s&=\sum_{\ell=0}^L s^\ell h_\ell,
 \qquad \Xi_s=\sum_{j=0}^{r-1}f^{s+1-j}\eta_j(s).
 \end{aligned}
\end{equation}
\end{lemma}

\begin{proof}
The Leibniz rule gives
$$
 (\partial_{x_i}u)v\omega
 =-u(\partial_{x_i}v)\omega
   +d\bigl(uv\,\iota_{\partial_{x_i}}\omega\bigr),
$$
where $\iota_{\partial_{x_i}}\omega=(-1)^{i-1}dx_1\wedge\cdots\wedge
\widehat{dx_i}\wedge\cdots\wedge dx_n$.  Apply this identity successively
to each summand of $P(s)f^{s+1}$, transferring every derivative from
$f^{s+1}$ to the coefficient of the differential operator.  The term
outside the total differential is then $f^{s+1}h_s\omega$; differentiating
these coefficients does not increase their degree bound $L$ in $s$.
Each term inside the total differential contains at most $r-1$
derivatives of $f^{s+1}$.  The rule
$\partial_{x_i}f^{s+1}=(s+1)(\partial_{x_i}f)f^s$ shows inductively
that such derivatives are sums of $f^{s+1-j}$ times holomorphic
coefficients polynomial in $s$, with $0\le j\le r-1$.
This gives the stated expression for $\Xi_s$.
\end{proof}

Here and below, a smooth singular $(n-1)$-chain in a Milnor fiber is a
finite complex linear combination of smooth maps of real simplices
$\Delta^{n-1}\to f^{-1}(t)$. A smooth family
of chains over $t\in U$ is a finite linear combination with complex coefficients (constant in $t$) of fibres $\Delta^{n-1}\to f^{-1}(t)$ of smooth maps $U\times \Delta^{n-1}\to f^{-1}(U)$.

\begin{lemma}[Absolute-period criterion]\label{lem:bs-period}
Let $f:(\C^n,x)\to(\C,0)$ be a holomorphic germ, and let $\gamma_t$,
$0<t\le\epsilon$, be a smooth family of $(n-1)$-cycles with
complex coefficients in its local Milnor fibers.  Suppose that their
supports tend uniformly to $x$ as $t\to 0$.  Put
$$
 I_h(t)=\int_{\gamma_t}\frac{h\omega}{df},
 \qquad h\in\mathcal O_{\C^n,x}.
$$
If $q>0$ and
\begin{equation}\label{eq:bs-all-numerators}
 I_h(t)=O_h(t^{q-1})\quad\hbox{for every $h$,\footnotemark}
 \qquad I_1(t)=t^{q-1}(c+o(1)),\quad c\ne0,
\end{equation}
\footnotetext{The subscript $h$ means that $|I_h(t)|\le C_h t^{q-1}$ for all
sufficiently small $t>0$, where the constant $C_h$ and the required
smallness of $t$ may depend on $h$.}
then $b_{f,x}(-q)=0$.  
\end{lemma}

\begin{proof}
Apply the identity of Lemma~\ref{lem:bs-leibniz} with
$b=b_{f,x}$, and choose representatives of its coefficients on a common
neighborhood of $x$.  Decrease $\epsilon$ so that the cycles $\gamma_t$
for $0<t\le\epsilon$ lie in this neighborhood.  Define
$$
 M_h(s)=\int_0^\epsilon t^sI_h(t)\,dt.
$$
For $h=1$ and the finitely many $h_\ell$ in that lemma, the first estimate
in~\eqref{eq:bs-all-numerators} makes $M_h$ holomorphic for
$\mathop{\rm Re}s>-q$.

For $0<a<\epsilon$, the family $\gamma_t$ over $a\le t\le\epsilon$
sweeps out an $n$-chain with boundary $\gamma_\epsilon-\gamma_a$.
Integrating \eqref{eq:bs-green-identity} over it and applying Stokes'
theorem gives
$$
 b(s)\int_a^\epsilon t^sI_1(t)\,dt
 =\int_a^\epsilon t^{s+1}I_{h_s}(t)\,dt
   +B_\epsilon(s)-B_a(s),
 \qquad B_t(s)=\int_{\gamma_t}\Xi_s.
$$
For each holomorphic $(n-1)$-form $\eta$, write
$df\wedge\eta=g\omega$.  Its restriction to the fiber is
$g\omega/df$, so its period is $I_g(t)=O_g(t^{q-1})$.
Applying this to the finitely many coefficients of the $\eta_j(s)$
shows that $B_a(s)\to0$ as $a\to0$ whenever
$\mathop{\rm Re}s>r-1-q$.  The fixed endpoint term $B_\epsilon(s)$ is
entire in $s$.  Using the expression for $h_s$ from
Lemma~\ref{lem:bs-leibniz}, we obtain
\begin{equation}\label{eq:bs-mellin-identity}
 b(s)M_1(s)=\sum_{\ell=0}^L s^\ell M_{h_\ell}(s+1)
                 +B_\epsilon(s).
\end{equation}
Initially this holds for sufficiently large $\mathop{\rm Re}s$.
Both sides are holomorphic for $\mathop{\rm Re}s>-q$, so the identity
holds throughout that half-plane.  Its right side is moreover
holomorphic for $\mathop{\rm Re}s>-q-1$.
The second estimate in \eqref{eq:bs-all-numerators} implies
$$
 \lim_{u\downarrow0}uM_1(-q+u)=c\ne 0.
$$
To see this, write $t^{1-q}I_1(t)=c+r(t)$, where $r(t)\to0$.
The constant term contributes $c\epsilon^u$, while
$u\int_0^\epsilon t^{u-1}r(t)\,dt\to0$ by splitting the integral at a
small fixed positive number.  Multiply \eqref{eq:bs-mellin-identity}
by $u$, put $s=-q+u$, and let $u\to 0$ from above.  The right side tends to
zero because it is holomorphic at $-q$, whereas the left side tends to
$b(-q)c$.  Since $c\ne0$, this proves $b(-q)=0$.
\end{proof}

\begin{rem}
For the Leray tubes used below, the first estimate in
\eqref{eq:bs-all-numerators} cannot fail while the second holds:
formula~\eqref{eq:tube-all-numerators} gives the stronger identity
$I_h(t)=t^{q-1}(h(x)c+o_h(1))$ for every $h$.
For arbitrary smooth families of cycles, however, the estimate for all
numerators is a separate hypothesis.  \end{rem}

We now apply the criterion to curve residues. Let $|I|=n-1$ and $\nu_i=qN_i$ for $i\in I$, so that $W_{I,q}$ is a
one-form on the $q$-effective cover of the curve stratum.

\begin{utver}[A curve residue period gives a root]
\label{prop:bs-residue-root}
Suppose that a finite complex linear combination $\gamma$ of compact
smooth loops on $\widehat E^\circ_{I,q}$ has projected support in
$\rho^{-1}(x)$ and satisfies
$$
                         \int_\gamma W_{I,q}\ne0.
$$
Then $b_{f,x}(-q)=0$.
\end{utver}

\begin{proof}
Take the full inverse image of $\gamma$ under
$\widehat E_I^\circ\to\widehat E_{I,q}^\circ$, with the inherited
orientations and coefficients.  It is a finite complex linear combination
of compact smooth loops, and its residue period is
$(m_I/d)\int_\gamma W_{I,q}\ne0$, where $d$ is the denominator of $q$
in lowest terms.  Its Leray tubes $\gamma_t$ form a smooth family of
cycles over $t>0$, and their supports tend uniformly to $x$.
Formula~\eqref{eq:tube-all-numerators} gives, for every holomorphic
germ $h$ at $x$,
\begin{equation}\label{eq:bs-leading-period}
 I_h(t)=t^{q-1}\bigl(h(x)c+o_h(1)\bigr),\qquad
 c=\pm\frac{(2\pi i)^{n-2}}{d}\int_\gamma W_{I,q}\ne0.
\end{equation}
These are the estimates in Lemma~\ref{lem:bs-period}, so that lemma proves
the assertion.
\end{proof}

\begin{sledst}\label{cor:bs-individual}
Let $C\subset\rho^{-1}(0)$ be a clean vertical equal-ratio curve.
If the contribution $B_{\mathrm{mot}}(C)$ in~\eqref{eq:res-Bmot}
is nonzero (in particular, if $B_{\mathrm{top}}(C)\ne0$), or $C$ has genus one
and all its defects are $1$, then $b_{f,0}(-q)=0$.  \end{sledst}

\begin{proof}
Theorem~\ref{thm:res-individual-RV-curve} and
Corollary~\ref{cor:poles-individual-motivic} supply a nonzero compact
residue period under either hypothesis.  Apply
Proposition~\ref{prop:bs-residue-root}.
\end{proof}    

\subsection{Proof of the root theorem}
\label{subsec:bs-horizontal}

We use two properties of the local Bernstein--Sato polynomial.  First, restriction of an
identity gives
\begin{equation}\label{eq:bs-nearby-divisibility}
                  b_{f,x}(s)\mid b_{f,0}(s)
                  \qquad\text{for every sufficiently nearby }x.
\end{equation}

Second, from \cite[Theorem~10.6]{KollarPairs}, the largest root of $b_{f,x}$ is $-\lct_x(f)$, where
$$
 \lct_x(f)=\min_{x\in\rho(E_i)}\frac{\nu_i}{N_i}.
$$
We moreover use the refinement from \cite[\S2.3 and Theorem~2.4]{NicaiseXu}:
if $C$ is a connected component of $E_I$ with $q={\nu_i}/{N_i},\, i\in I$, and
\begin{equation}\label{eq:bs-NX-hypotheses}
 C\cap\rho^{-1}(x)\ne\varnothing,\qquad
 C\cap\rho^{-1}(x)\cap E_j=\varnothing\quad(j\notin I),
\end{equation}
then $q=\lct_x(f)$.

\begin{utver}\label{prop:bs-horizontal}
If a horizontal component $C$ of the curve $E_I$, with
$q=\nu_i/N_i$ for all $i\in I$, meets $\rho^{-1}(0)$, then
$b_{f,0}(-q)=0$.
\end{utver}

\begin{proof}
For $j\notin I$, the intersection $C\cap E_j$ is
finite by simple normal crossings.  For $x\in\rho(C)$ outside their
finitely many images, the conditions
\eqref{eq:bs-NX-hypotheses} hold.  Hence
$q=\lct_x(f)$ and $b_{f,x}(-q)=0$.  Such points occur arbitrarily
close to $0$; \eqref{eq:bs-nearby-divisibility} transfers their root
to $b_{f,0}$.
\end{proof}

\begin{proof}[Proof of Theorem~\ref{thm:bs-high-order}]
Theorem~\ref{thm:poles-global-RV} supplies a horizontal $q$-curve
meeting $\rho^{-1}(0)$ or a clean vertical $q$-curve satisfying the hypotheses of Corollary~\ref{cor:bs-individual}.  In the
horizontal case apply Proposition~\ref{prop:bs-horizontal}; in the
vertical case apply that corollary.  \end{proof}

\end{document}